\documentclass[12pt,leqno]{article}

\usepackage{amssymb}
\usepackage{amsmath}
\usepackage{color}
\usepackage{url}

\newtheorem{thm}{Theorem}[section]
\newtheorem{defi}{Definition}[section]

\newtheorem{lem}{Lemma}[section]
\newtheorem{rem}{Remark}[section]
\newtheorem{cor}{Corollary}[section]

\newcommand{\comment}[1]{}

\newcommand{\de}{:=}
\renewcommand{\a}{\alpha}
\newcommand{\Aa}{\mbox{$\mathfrak A$}}
\newcommand{\Bb}{\mbox{$\mathfrak B$}}
\newcommand{\Cc}{\mbox{$\mathfrak C$}}
\newcommand{\Ff}{\mbox{$\mathfrak F$}}
\newcommand{\Dd}{\mbox{$\mathfrak D$}}
\newcommand{\Mm}{\mbox{$\mathfrak M$}}

\newcommand{\Rd}{\mbox{$\mathfrak R\mathfrak d$}}
\newcommand{\Rrd}{\mbox{\sf Rd}}
\newcommand{\Nr}{\mbox{$\mathfrak N\mathfrak r$}}
\newcommand{\Th}{\mbox{\sf Th}}
\newcommand{\CA}{\mbox{\sf CA}}
\newcommand{\Pow}{\mbox{$\mathcal P$}}
\newcommand{\Ll}{\mbox{$\mathcal L$}}

\newcommand{\Ws}{\mbox{\sf Ws}}
\newcommand{\RCA}{\mbox{\sf RCA}}
\newcommand{\Cs}{\mbox{\sf Cs}}

\newcommand{\iCs}{\mbox{$\Iso\,{}_\infty\!\Cs$}}
\newcommand{\Di}{\mbox{\sf Di}}
\newcommand{\Sy}{\mbox{\sf Sy}}
\newcommand{\Ret}{\mbox{\sf Ind}}
\newcommand{\Edc}{\mbox{\sf Edc}}
\newcommand{\Dc}{\mbox{\sf Dc}}
\newcommand{\Lf}{\mbox{\sf Lf}}

\newcommand{\PEA}{\mbox{\sf PEA}}

\newcommand{\QPEA}{\mbox{\sf QPEA}}
\newcommand{\FPEA}{\mbox{\sf FPEA}}
\newcommand{\K}{\mbox{\sf K}}
\newcommand{\LL}{\mbox{\sf L}}
\newcommand{\V}{\mbox{\sf V}}
\newcommand{\Eq}{\mbox{\sf Eq}}
\newcommand{\QED}{\hfill\mbox{$\Box$}\bigskip}
\newcommand{\s}{\mbox{\sf s}}
\newcommand{\p}{\mbox{\sf p}}
\newcommand{\Dom}{\mbox{\sf Dom}}
\newcommand{\Rg}{\mbox{\sf Rg}}
\newcommand{\ind}{\mbox{\sf ind}}
\newcommand{\Sub}{\mbox{\bf S}}
\newcommand{\Iso}{\mbox{\bf I}}

\newcommand{\Var}{\mbox{\bf Var}}

\begin{document}

\title{How many varieties of cylindric algebras}
\author{Andr\'eka, H.\ and N\'emeti, I.}%
\date{}
\maketitle
\begin{abstract}
Cylindric algebras, or concept algebras in another name, form an
interface between algebra, geometry and logic; they were invented by
Alfred Tarski around 1947. We prove that there are $2^\alpha$ many
varieties of geometric (i.e., representable) $\alpha$-dimensional
cylindric algebras, this means that $2^\alpha$ properties of
definable relations of (possibly infinitary) models of first order
logic theories can be expressed by formula schemes using $\alpha$
variables, where $\alpha$ is infinite. This solves Problem 4.2 in
the 1985 Henkin-Monk-Tarski monograph \cite{HMTII}, the problem is
restated in \cite{NAPAL87, AMonkN91}. For solving this problem, we
had to devise a new kind of construction, which we then use to solve
Problem 2.13 of the 1971 Henkin-Monk-Tarski monograph \cite{HMTI}
which concerns the structural description of geometric cylindric
algebras. There are fewer varieties generated by locally finite
dimensional cylindric algebras, and we get a characterization of
these among all the $2^\alpha$ varieties. As a by-product, we get a
simple, natural recursive enumeration of all the equations true of
geometric cylindric algebras, and this can serve as a solution to
Problem 4.1 of the 1985 Henkin-Monk-Tarski monograph. All this has
logical content and implications concerning ordinary first order
logic with a countable number of variables.
\end{abstract}

\section{Introduction}
Cylindric algebras (or concept algebras in another name) are an
algebraic form of first order logic, analogous to Boolean algebras
which are an algebraic form of propositional logic. Cylindric
algebras were created by Alfred Tarski around 1947, they are Boolean
algebras endowed with complemented closure operators---one for each
quantifier $\exists v_i$---and constants $v_i=v_j$ for representing
equality. Set cylindric algebras analogous to set Boolean algebras
are called representable cylindric algebras, these latter are
equationally definable (i.e., form a variety), this was proved by
Tarski in 1952. However, unlike the propositional case, not all
cylindric algebras are representable, and this cannot be repaired
easily since J. Donald Monk proved in 1969 that the representable
algebras are not finite schema axiomatizable. This gap between
``abstract" and representable cylindric algebras is in the center of
algebraic logic studies, it is a source of insights and problems.

A set cylindric algebra is the algebra of all definable relations of
a model. Representable cylindric algebras are the same, just for
classes of models, i.e., for theories, in place of single models. We
often call them concept algebras, since they are natural algebras of
concepts of a theory. (The name ``cylindric algebras" on the other
hand refers to the geometrical meaning of these algebras.)
Properties of a model, or of a theory, can be read off its concept
algebra. E.g., it can be expressed by an equation whether dense
linear order can be defined in the model. Equations of concept
algebras talk about properties of all definable relations in a
model, as opposed to talking about the primitive relations only. An
equation in the algebraic language corresponds to a formula schema
of first order logic where the formula variables range over formulas
with free variables of the language, this way they ``talk about"
definable relations of models. In this context, the problem asking
about the number of varieties of representable cylindric algebras
asks how many properties of definable relations we can express by
first order logic
schemata. 
In this connection, the present work is not unrelated to Shelah's
classification theory. (For some connections, see \cite{SagiSziraki,
ASSS13}.)

Let $\a$ denote the number of variables we have in the first order
language the algebraic versions of which we investigate. In this
paper we deal with infinite $\a$ only, let $\RCA_\a$ denote the
class of representable cylindric algebras of this language. Then we
have $\a$ many closure operations called cylindrifications, and
$\a\times\a$ many constants called diagonal constants apart from the
Boolean operations, so our equational language has size $\a$ which
means that there are at most $2^{\a}$ many different equational
theories in the algebraic language. It was known that $\RCA_\a$ has
at least continuum many subvarieties (\cite[Thm.4.1.24]{HMTII}) and
\cite[Problem 4.2]{HMTII} asks whether there are $2^{\a}$ many for
uncountable $\a$. There are theorems pointing to the answer being
continuum (i.e., $2^{\omega}$, where $\omega$ is the least infinite
ordinal) and there are theorems pointing to the answer being
$2^{\a}$. For example, it is proved in \cite{NAPAL87} that in
concept algebras of finite models we can only express the size of
models nothing else, in the algebraic language the theorem says that
the concept algebras of models of size $n$ generate a variety which
is an atom of the lattice of varieties of $\RCA_\a$. On the other
hand, we can say lots of things about definable relations in
infinite models. E.g., the equational theory of the concept algebra
of a single infinite model in which we have no primitive relations
(i.e., in which all the definable relations are the ones we can
define by using the equality only) is not recursively enumerable,
see \cite[Thm.1(i)]{NAPAL87}. We prove in this paper that there are
$2^{\a}$ many subvarieties of $\RCA_\a$, but we also show that this
is only part of the answer since there is a large subclass of the
distinguished representable cylindric algebras which has only
continuum many subvarieties.

From the logical point of view, at first, the subject of the above
discussed problem might look ``esoteric", since it concerns the
difference between using countably or uncountably infinite number of
variables in the logical languages having only finite formulas.
However, our solution has impact on the countable case, too. We
mentioned already that the gap between abstract and representable
cylindric algebras is a source of insights and problems. One of the
insights has been that this gap is concerned with the number of
``free" extra variables $v_i$ we have for a formula or proof in the
logical language. The crucial thing is not that each first order
formula uses only finitely many variables, but instead that for each
finite set of formulas there is always at least one variable that
they do not use. Let us call a cylindric algebra
dimension-complemented if to each finite subset there is a
cylindrification all elements of this finite subset are fixed-points
of; one can prove then that all dimension-complemented cylindric
algebras are representable. The strongest and most beautiful form of
this insight is Leon Henkin's theorem saying that an abstract
cylindric algebra is representable if and only if it can be embedded
to one having infinitely many more cylindric operations so that each
of the elements of the original algebra are fixed-points of all the
new operations. Efforts were made to find structural properties of
abstract cylindric algebras, similar to being dimension-complemented
mentioned above, which refer to only the original cylindric algebra
without comparing it to others. Theorem 2.6.50 in \cite{HMTI}
summarizes how far they got in this respect, and \cite[Problem
2.13]{HMTI} asks if their last description (which we call here
endo-dimension-complemented) captures being representable or not.

Our full answer to \cite[Problem 4.2]{HMTII} shows that when proving
that there are many varieties of $\RCA_\a$, we had to use unusual
features of representable cylindric algebras, we had to devise a new
kind of construction. This contributes to understanding the
structural properties of representable cylindric algebras. Namely,
our construction shows that ``endo-dimension-complemented" is not
the final answer in \cite[Thm.2.6.50]{HMTI}, since our new
constructions are representable yet not endo-dimension-complemented,
solving Problem 2.13 in the negative. However, here, too, something
positive can be added as a result of this solution. A new property
extending ``endo-dimension-complemented" toward representability
emerges by an analysis of our construction. We call an abstract
cylindric algebra inductive iff from the fact that an equation holds
for all elements which are fixed-points of a given cylindric
operation, we can infer that this equation holds for the whole
algebra, provided that the given cylindrification does not occur in
the equation. In a way, this property also talks about extra free
variables. We prove that each endo-dimension-complemented algebra is
inductive and each inductive algebra is representable, but none of
the reverse implications holds. This way we get a new
representability theorem extending the chain of classes in
\cite[Thm.2.6.50]{HMTI}, and we get a new insight about
representable algebras.

A cylindric algebra is called locally finite dimensional if each of
its elements is the fixed-point of all but finitely many
cylindrifications, $\Lf_\a$ denotes their class. These algebras
correspond to theories of first order logic in which all primitive
relations have finite rank. The rest of the representable algebras
correspond to theories of first order logic when the primitive
relations may have arbitrary ranks, let us call this logic
infinitary first order logic to distinguish it from the former which
we call finitary or ordinary (following \cite{HMTII}). It was known
that the infinitary and finitary first order logics have the same
valid formulas, in algebraic form this theorem says that $\Lf_\a$
generates the variety $\RCA_\a$. However, it was not known whether
the two logics have the same theories, i.e., whether all
subvarieties of $\RCA_\a$ are generated as varieties by their
locally finite-dimensional members. In this paper we give a negative
answer to this: there are subvarieties of $\RCA_\a$ not generated by
their locally finite-dimensional members, even in the $\a=\omega$
case, so the theories of finitary and infinitary first order logics
are different.  The notion of inductivity is the key notion to
finding this answer.
The property of a cylindric algebra being inductive can easily be
converted to a new logical rule, which we call inductive rule. We
prove that this rule is admissible for ordinary first order logic
but it is not admissible for infinitary first order logic. In
particular, we prove that a theory is one of ordinary first order
logic iff it is closed under the inductive rule.

The above theorem readily yields a simple, natural recursive
enumeration for all equations valid in $\RCA_\a$. This gives a
possible solution to \cite[Problem 4.1]{HMTII}, this problem is
restated in numerous other places, e.g., in \cite[Problem 19,
p.735]{AMonkN91}, or in logical form as \cite[Problem 2.9]{NAPAL87}.

The structure of the paper is the following. In the rest of the
introduction we recall the definitions, notation, necessary
background needed in the paper. Section~\ref{manyvar-sec} contains
the proof that there are $2^{\a}$ varieties of $\RCA_\a$,
section~\ref{ca-sec} contains the proof of an analogous statement
for abstract cylindric algebras. Section~\ref{counter-sec} concerns
with subclasses where there are only continuum many varieties and
with the results we can get from the proof contained in
section~\ref{manyvar-sec}. We introduce the notions of symmetric,
endo-dimension-complemented (endo-dc in short), and inductive
algebras. Section~\ref{sym-ssec} investigates the relationship
between symmetric and endo-dc algebras leading to the solution of
\cite[Problem 2.13]{HMTI}. Section~\ref{poly-ssec} reveals
connections with the polyadic substitution operations.
Section~\ref{ind-ssec} concerns inductive algebras. It places this
notion between endo-dc and symmetric and shows representability of
inductive algebras by proving that they are exactly the algebras
equationally indistinguishable from a member of $\Lf_\a$. As a
corollary we get a $\Delta_2$-formula separating $\Lf_\a$ and
$\RCA_\a$. Section~\ref{rcaeq-ssec} contains a simple recursive
enumeration of the equational theory of $\RCA_\a$ and the
characterization of varieties generated by subclasses of $\Lf_\a$.
Finally, section~\ref{fewvar-ssec} contains a proof that locally
finite dimensional set algebras with infinite bases have indeed
continuum many subvarieties.

\subsection{Background}\label{back-ssec}

Where not specified otherwise, we use the notation of \cite{HMTI,
HMTII}, but we try to be self-contained.

Let $\a$ be any set. An algebra $\langle A,+,-,c_i, d_{ij} :
i,j\in\a\rangle$ is a \emph{cylindric algebra} of dimension $\a$ if
$\langle A,+,-\rangle$ is a Boolean algebra and for all distinct
$i,j,k\in\a$ the following hold. The operations $c_i$ are unary,
they are commuting complemented closure operators, i.e., for all
$x,y\in A$ we have $c_ic_jx=c_jc_ix$, $x\le c_ix=c_ic_ix$,
$c_i(x+y)=c_ix+c_iy$, $c_i-c_ix=-c_ix$.  The operations $d_{ij}$ are
nullary, i.e., they are constants and they satisfy the following
equations: $d_{ii}=1$, $d_{ij}=d_{ji}$, $c_j(d_{ij}\cdot d_{jk})=
d_{ik}$, $c_id_{ij}=1$, $d_{ij}\cdot c_i(d_{ij}\cdot x)=d_{ij}\cdot
x$. In the above, $\cdot$ is the Boolean intersection defined from
$+,-$ the usual way. The extra-Boolean operations $c_i$ and $d_{ij}$
are called \emph{cyindrifications} and \emph{diagonal constants},
respectively. The schemata of equations $(C_0) - (C_7)$ in
\cite[p.161]{HMTI} express the above. $\CA_\a$ denotes the class of
all cylindric algebras.

The cylindric \emph{set algebras} of dimension $\a$ are Boolean set
algebras of $\a$-dimensional spaces, where the extra-Boolean
operations have natural geometric interpretations. Let $U$ be any
set, the points of the $\a$-dimensional space ${}^\a U$ over $U$ are
the $U$-termed $\a$-sequences, i.e., the set of all functions from
$\a$ to $U$. We use the terms sequence and function to mean the same
thing. If $s$ is any sequence, then $\Dom(s), \Rg(s)$ denote its
domain and range, $s_i$ denotes $s(i)$, and $s(i/u)$ denotes the
sequence we get from $s$ by changing its value at $i$ to be $u$ if
$i\in\Dom(s)$ (and $s(i/u)$ denotes $s$ if $i\notin\Dom(s)$).
For $V\subseteq {}^\a U$ the \emph{full cylindric set algebra with
unit $V$} is
\begin{description}
\item{}
$\langle \Pow(V), \cup, - , C_i^V, D_{ij}^V : i,j\in\a\rangle$\quad
where $\Pow(V)$ is the powerset of $V$ and
\item{}
$C_i^VX \de \{ s(i/u)\in V : s\in X, u\in U\}$,\quad $D_{ij}^V \de
\{ s\in V : s_i=s_j\}$.
\end{description}
Algebras isomorphic to subalgebras of full cylindric set algebras
with units as disjoint unions of $\a$-dimensional spaces are called
\emph{geometric}, or \emph{representable}, $\RCA_\a$ denotes their
class. (For technical reasons, when $\a$ is a one-element set,
$\RCA_\a$ is defined as the class of subdirect products of these.)
Cylindric set algebras with units of form ${}^\a U$ are called
simply \emph{cylindric set algebras}, $U$ is called their \emph{base
set}, and their class is denoted as $\Cs_\a$.

Geometric cylindric algebras have natural \emph{logical meaning},
too. Let us have a first order equality language $\Ll$ with
variables as $\a$ (or $v_i$ with $i\in\a$), and with the logical
connectives $\lor, \neg, \exists v_i, v_i=v_j$ for $i\in\a$, and let
$\Mm$ be a relational structure with universe $U$ and with at most
$\a$-place primitive relations. The points of the $\a$-dimensional
space over $U$ are also evaluations of variables. For any formula
$\varphi$ in this language let $\varphi^{\Mm}$ denote the set of
evaluations of variables under which $\varphi$ is true in $\Mm$.
Then $(\varphi\lor\psi)^{\Mm} = \varphi^{\Mm}\cup\psi^{\Mm},
(\neg\varphi)^{\Mm}={}^\a U - \varphi^{\Mm}, (\exists
v_i\varphi^{\Mm})=C_i\varphi^{\Mm}$ and $(v_i=v_j)^{\Mm}=D_{ij}$
(with $V={}^\a U$). Thus
\[ \mbox{Ca}^{\Mm} \de \{\varphi^{\Mm} : \varphi\in\Ll\} \]
is a subuniverse of the full cylindric set algebra with unit ${}^\a
U$. We call the subalgebra with this universe the \emph{concept
algebra} of $\Mm$, its universe is the set of all definable
relations over $\Mm$. Concept algebras for theories can be defined
analogously, the unit of such a concept algebra is a disjoint union
of sets of form ${}^\a U$, so it is in $\RCA_\a$. Two theories are
definitionally equivalent iff their concept algebras are isomorphic,
and homomorphisms from one concept algebra to another correspond
exactly to the interpretations between the two theories. Thus, the
category $\RCA_\a$ with homomorphisms corresponds to the category of
all theories and interpretations between them (see
\cite[sec.4.3]{HMTII}). This feature of concept algebras comes handy
in applications of logic, e.g., in physics. For this kind of
applications see, e.g., \cite{ANComp, BaHal15, LefSzek, MadDis,
MadSzek, Weatherall}.

Equations holding for concept algebras have several logical
interpretations, see e.g., \cite[sec.4.3]{HMTII}. The most direct of
these is where we interpret the algebraic variables occurring in the
equation as schemes for formulas in the corresponding logical form.
This way we get the logic of formula schemata, for definitions see
e.g., \cite{NAPAL87}, \cite[sec.3.7]{Rybakov}. This schema logic
talks directly about the definable relations in a model.
Alternatively, we can interpret the algebraic variables as primitive
$\a$-place relation symbols (full restricted first order logic in
\cite{HMTII}), or as primitive relation symbols of unspecified but
finite arity (type-free logic in \cite{HMTII}). Via these logical
interpretations, results about varieties of $\RCA_\a$ imply
corresponding results about these logics. We deal with the algebraic
aspects in this paper, the logical consequences are dealt with in a
separate paper.

A class of algebras that can be axiomatized/defined by a set of
equations is called a \emph{variety}. Let $\K$ be a class of
algebras of the same similarity class. Then $\Eq(\K)$ denotes the
set of all equations (using a countable set of prespecified
algebraic variables) valid in all algebras in $\K$, and $\K$ is a
variety iff $\K$ consists of all algebras in which $\Eq(\K)$ is
true. The variety generated by $\K$ is the least variety containing
$\K$, this is the class of all algebras in which $\Eq(\K)$ holds.
Varieties are one of the main subjects of universal algebra,
Birkhoff's theorem says, e.g., that all members of the variety
generated by $\K$ can be obtained from members of $\K$ as
homomorphic images of subalgebras of products of members of $\K$.
The class of algebras isomorphic to members of $\K$ is denoted by
$\Iso\K$, the class of subalgebras of members of $\K$ is denoted by
$\Sub\K$.

The following are the main facts, in connection with the present
paper, known about varieties of cylindric algebras, these are
contained in \cite{HMTI, HMTII, HMTAN} if not specified otherwise.
$\RCA_\a\subseteq\CA_\a$ is a variety, it is not finitely
axiomatizable iff $|\a|>2$, where if $\a$ is any set, $|\a|$ denotes
its cardinality. J. D. Monk \cite{Monk69, Monk70} characterized the
lattice of all subvarieties of $\RCA_\a$ for $|\a|=1$, and a similar
characterization for $|\a|=2$ is contained in \cite{BezhCAIII}. Let
$\a$ be infinite. Some of the distinguished subvarieties of
$\RCA_\a$ are $\iCs_\a$ and ${}_n\RCA_\a$ for finite $n$. $\iCs_\a$
is the class of all algebras isomorphic to a cylindric set algebra
with unit of form ${}^\a U$ for infinite $U$, and ${}_n\RCA_\a$ is
the class of all algebras isomorphic to cylindric set algebras with
unit as disjoint union of sets of form ${}^\a U$ where $|U|=n$. It
is proved in \cite{NAPAL87} that ${}_n\RCA_\a$ are atoms in the
lattice of subvarieties of $\RCA_\a$, but $\iCs_\a$ is not an atom.
The structure of subvarieties is interesting and is investigated for
other kinds of algebras related to logic as well, see e.g.,
\cite{Blok, Blok2, Jonsson, Goldblatt, JiRo, AGiNJSL94}.

Monographs and books on these algebras and their logical
applications include \cite{Hal62, HMTI, Craig, HMTAN, HMTII,
AMonkN91, Marx, Rybakov, HHbook, GKWZ, Madduxbook, AFerNCAIII,
Givantbook}.

\goodbreak

\section{Number of subvarieties of $\RCA_\a$}\label{manyvar-sec}

Let $\a$ be any infinite ordinal, throughout the rest of the paper.
(We assume that $\a$ is an ordinal, and not just any set, for
convenience, this way $\a$ is ordered by the elementhood relation.)
It is proved in \cite[4.1.24]{HMTII} that there are at least
$2^\omega$ many subvarieties of $\RCA_\a$. Since in the language of
$\CA_\a$ there are $|\a|$ many equations, there can be at most
$2^{|\a|}$ many subvarieties of any $\K\subseteq \CA_\a$. Problem
4.2 in \cite{HMTII}  asks if there are $2^{|\a|}$ many subvarieties
of $\RCA_\a$ and of $\Iso\Cs_\a\subseteq\RCA_\a$, if $\a$ is
uncountable. The problem is restated in \cite[Problem 41,
p.738]{AMonkN91}. In this section, we prove that indeed there are
maximum number of subvarieties of $\RCA_\a$ as well as of
$\iCs_\a\subseteq\Iso\Cs_\a$. Note that both $\RCA_\a$ and $\iCs_\a$
are varieties but $\Iso\Cs_\a$ is not.

\begin{thm}\label{manyvar-thm}\text{\rm (Solution of \cite[Problem
4.2]{HMTII})} Let $\a$ be infinite. There are $2^{|\a|}$ many
subvarieties of $\RCA_\a$ as well as of $\iCs_\a$.
\end{thm}

\noindent The proof of Theorem~\ref{manyvar-thm} is contained in
subsections~\ref{const-sec}-\ref{check-sec}. The idea of the proof
is the following. We exhibit a set of equations $E$ of cardinality
$|\a|$ which is independent in the sense that no element of $E$
follows from the rest of the equations in $E$. All these equations
will be variants of a single equation $e$ such that we rename the
indices of the operations occurring in $e$. We will show
independence of $E$ by constructing one algebra $\Aa\in\iCs_\a$ in
which $e$ fails but the rest of the equations in $E$ hold. Then the
algebras in which we rename the operations $c_i$ and $d_{ij}$
according to appropriate permutations of $\a$ will show independence
of the whole set $E$. This will show that all the subvarieties of
$\iCs_\a$ specified by the $2^{|\a|}$ many subsets of $E$ are
distinct. We begin with constructing the ``witness" algebra $\Aa$,
because it will give intuition for writing up the ``master" equation
$e$.

\subsection{Construction of the witness algebra $\Aa$}\label{const-sec}
Let $\langle V_i : i\in\a\rangle$ be a system of sets such that
\[\mbox{$V_0=V_1=V_2$ is the set of rational numbers,}\]
and all the other $V_i$'s are pairwise disjoint two-element sets
disjoint from $V_0$, too (i.e., $V_i\cap V_j=\emptyset$ for $2\le
i<j<\a$). Let $U$ be the union of these sets, i.e.,
\[ U\de\bigcup\{ V_i : i\in\a\}, \]
let $p$ be an $\a$-sequence such that $p_i\in V_i$ for all $i\in\a$
and let $V$ be the set of $U$-termed $\a$-sequences that deviate
from $p$ only at finitely many places, i.e.,
\[ V\de {}^\a U^{(p)} = \{ s\in {}^\alpha U : |\{ i\in\a :s_i\ne p_i\}|<\omega\}. \]
This $V$ will be the unit of our algebra. Since $V_0$ is the set of
rational numbers, we will use the usual operations and ordering $<$
between rational numbers. Our algebra is generated by a single
element, namely by
\begin{description}
\item{}
$g\de\{s\in V : s_0<s_1<s_2\mbox{ and } s_i\in V_i\mbox{ for all }i\in\a\mbox{ and }\\
s_1=(s_0+s_2)\slash 2\mbox{\quad if }|\{ s_i : s_i\ne p_i, i>2\}| \mbox{ is even,}\\
s_1\ne(s_0+s_2)\slash 2\mbox{\quad if }|\{ s_i : s_i\ne p_i,
i>2\}|\mbox{ is odd \ }\}.$
\end{description}
Let $\Aa$ denote the subalgebra of the full set algebra with unit
$V$ that is generated by $g$.

Sets of form ${}^\a U^{(p)}$ for some set $U$ and $p\in{}^\a U$ are
called \emph{weak spaces} and algebras with unit a weak space are
called \emph{weak set algebras}, their class is denoted by $\Ws_\a$.
It is proved in \cite[3.1.102]{HMTII} that
$\Ws_\a\subseteq\Iso\Cs_\a$, thus our above constructed algebra
$\Aa$ is in $\iCs_\a$ since the fact that $U$ is infinite is
reflected by an equation holding in $\Aa$ (see,
\cite[2.4.61]{HMTI}). (We note that we could have used for our
witness algebra the subalgebra of the full cylindric set algebra
with unit ${}^\a U$ generated by $g$, the proofs would be only
slightly more complicated.) The set $U$ is called the \emph{base
set} of $\Aa$.

\subsection{Describing the elements of $\Aa$}\label{elem-sec}
The idea behind the construction of the witness algebra $\Aa$,
defined in the previous section, is that we put some information in
$g$ at the indices $0,1,2$ which information cannot be transferred
by the cylindric operations to higher indices $i,j,k\in\a$. Then we
express lack of this information by an equation $e$. If we succeed
with realizing this idea, then $e$ would fail in $\Aa$ at $0,1,2$
while $e$ would be valid in $\Aa$ at higher indices. We now turn to
elaborating this idea.

Let $R\subseteq {}^nU$ be an $n$-place relation on $U$, where $n$ is
any ordinal. We say that $X$ is a \emph{sensitive cut} of $R$ if
$c_iX=c_i(R-X)$ for all $i<n$. (Here, $c_i$ denotes cylindrification
of the full $\Cs_n$ with base set $U$, i.e., $c_iX=\{ s(i/u)\in
{}^nU : s\in X\}$.) Thus, as soon as we apply a cylindrification to
$X$, the information on how $X$ cuts $R$ into two parts is lost.
This technique is widely applicable, cf., e.g., \cite{ACMNS, Sayed,
Sereny, ST}.
We are going to show that our generator $g$ is a sensitive cut of
\[ T \de \{ s\in V : s_0<s_1<s_2\mbox{ and } s_i\in V_i\mbox{ for all
}i\in\a\} .\] Intuitively, the proof will be a kind of ``flip-flop"
play between the two independent conditions $s_1=(s_0+s_2)/2$ and
$|\{ s_i : s_i\ne p_i, i>2\}|$ being even in the definition of $g$.
Indeed, let $i\in\a$ and $s\in T$. We show that $s(i/u)\in g$ while
$s(i/v)\in T-g$ for some $u,v\in U$. This will show that
$c_ig=c_i(T-g)$. By $s\in T$ we have $s_0<s_1<s_2$. Let $\Sigma=|\{
s_i : s_i\ne p_i, i>2\}|$, $u=(s_0+s_2)\slash 2$, and let $v$ be
such that $v\ne u$, $s_0<v<s_2$. Assume first that $i=1$. Now, if
$\Sigma$ is even then $s(1/u)\in g$, $s(1/v)\in T-g$ and if $\Sigma$
is odd then $s(1/u)\in T-g$, $s(1/v)\in g$. For $i=0$ choose
$u=2s_1-s_2$ and $v<u$, for $i=2$ choose $u=2s_1-s_0$ and $v>u$, for
these choices the same is true as in the case of $i=1$. Assume $i>2$
and let $v$ be the element of $V_i$ distinct from $s_i$. Assume
$s_1=(s_0+s_2)/2$. Then  $s\in g$, $s(i/v)\in T-g$ if $\Sigma$ is
even, and otherwise $s\in T-g$, $s(i/v)\in g$. Assume
$s_1\ne(s_0+s_2)/2$. Then just the other way round, namely $s\in
T-g$, $s(i/v)\in g$ if $\Sigma$ is even, and otherwise $s\in g$,
$s(i/v)\in T-g$. We have shown that
\begin{equation}\label{cut-eq} c_ig = c_i(T-g) = c_iT\quad\mbox{ for all
}i\in\a .
\end{equation}
Next we show that this last property implies that the elements of
$A$ are those that are generated by $T$ and perhaps one of $g$,
$T-g$ added:
\begin{lem}\label{gen-lem} Let $\Bb$ be the weak set algebra of dimension $\a$ with unit $V$
and generated by $T$. Then
\begin{equation*}\label{algebra-eq}
A = \{ x+h :  x\in B\mbox{ and }h\in\{ 0,g,T-g\} \} .\end{equation*}
\end{lem}
\noindent {\bf Proof.} Each element of the form $x+h$ is generated
by $g$, since $T=c_0g\cdot c_2g$. To finish the proof, we are going
to show that $A$ is closed under the cylindric operations $c_i$ and
$d_{ij}$ as well as under the Boolean operations $+,-$. Let
$i,j\in\a$. Now, $A$ is closed under $c_i$ by \eqref{cut-eq} since
$c_i(x+h) = c_ix + c_ih\in B$ by $x\in B$ and $c_ih\in\{
0,c_iT\}\subseteq B$. Also, $d_{ij}\in A$ by $d_{ij}\in B$. The set
$A$ is closed under Boolean addition $+$ by its definition and by
$g+(T-g)=T\in B$.

To see that $A$ is closed under Boolean complementation, first we
show that $T$ is an atom in $B$. We will use the following property
of $<$ later on, too:
\begin{description}
\item{(*)}\label{ord-eq}
\quad Assume that $a_1<a_2<\dots<a_n$ and $b_1<b_2<\dots<b_n$. There
is an automorphism $\pi$ of $\langle V_0,<\rangle$  mapping
$a_1,\dots,a_n$ to $b_1,\dots,b_n$, respectively. If $a_1=b_1$ and
$a_n=b_n$ then $\pi$ can be chosen such that it is the identity on
elements smaller than $a_1$ or bigger than $a_n$.
\end{description}
To prove (*), for $a<b$ let $[a,b]=\{ x\in V_0 : a\le x\le b\}$
denote the closed interval between $a$ and $b$. Let $u,v\in V_0$ be
such that $u<a_0,\ u<b_0$ and $v>a_n,\ v>b_n$ and define $a_0\de
b_0\de u$ and $a_{n+1}\de b_{n+1}\de v$. For $k\le n$ and $x\in
[a_k,a_{k+1}]$ let $\pi_k(x)\de
(x-a_k)\cdot(b_{k+1}-b_k)\slash(a_{k+1}-a_k)+b_k$. Then $\pi_k$ is
an isomorphism between $\langle[a_k,a_{k+1}],<\rangle$ and
$\langle[b_k,b_{k+1}],<\rangle$, thus their union $\sigma\de
\pi_0\cup\dots\cup\pi_n$ is an automorphism of
$\langle[a_0,a_{n+1}],<\rangle$ taking $a_0,\dots,a_{n+1}$ to
$b_0,\dots,b_{n+1}$, respectively. We now can choose $\pi$ to be the
identity outside $[a_0,a_{n+1}]$ and $\sigma$ on the interval. This
proves (*).

Returning to showing that $T$ is an atom, let $s,z\in T$ be
arbitrary. Then $s_0<s_1<s_2$ and $z_0<z_1<z_2$ and $s_i,z_i\in V_i$
for all $i\in\a$ by the definition of $T$. Let $\pi$ be a
permutation of $U$ which takes $s_i$ to $z_i$ for all $i\in\a$, is a
permutation of $V_i$ for all $i\in\a$, and is an automorphism of
$\langle V_0,<\rangle$. By (*), there is such a $\pi$. Then clearly,
$\pi$ takes $s$ to $z$ in the sense that $z=\pi(s)\de \langle
\pi(s_k) : k\in\a\rangle$ while leaving $T$ as well as $V$ fixed,
i.e., $T=\pi(T)\de\{ \pi(q) : q\in T\}$ and $V=\pi(V)\de\{\pi(q) :
q\in V\}$. This implies that $\pi(x)=x$ for all $x\in B$ (see, e.g.,
\cite[3.1.36]{HMTII}), so $s\in x$ implies $z\in x$ for all $x\in
B$. Since $s,z\in T$ were chosen arbitrarily, this shows that $T$ is
an atom in $B$. We are ready to show that $A$ is closed under
complementation. Let $x\in B$ and $h\in\{ 0,g,T-g\}$. If $T$ is
disjoint from $x$ then $V-(x+h) = (V-x)+(T-h)\in A$, and if $T\le x$
then $x=x+h$ and we are done with proving Lemma~\ref{gen-lem}. \QED

\subsection{The set $E$ of independent equations}\label{equations-sec}
The ``master equation"  $e$ will express about an element that it is
not similar to our generator $g$. Namely, it will say about an
element $x$ that either it is not a sensitive cut of its closure
$c_0x\cdot c_2x$ (this is $T$ in the case of $g$), or else this
closure is not like $T$ in the sense that the first two coordinates
of $c_2x$ form a strict linear order $<_x$ and the ternary beginning
of the closure is $\{\langle u,v,w\rangle : u<_x v<_x w\}$. An
equation can talk about finitely many indices only, our equation
will concern the first three indices $0,1,2$.

We begin writing up the equation $e$.
First we write up a term we will use in checking that $x$ is not a
sensitive cut of its closure $z \de c_0x\cdot c_2x$ (precise
statements about the meanings of the terms below can be found in the
proof of Lemma~\ref{express-lem} that is to be stated soon).
\begin{equation*}
\beta(x)\de c_0x\oplus c_0(z-x) + c_1x\oplus c_1(z-x) + c_2x\oplus
c_2(z-x),
\end{equation*}
where $\oplus$ denotes Boolean symmetric difference, i.e., $x\oplus
y\de (x\cdot y)+(-x\cdot -y)$. In writing up the rest of the terms,
it will be convenient to use the following notation. It concerns
rearranging sequences in set algebras, see
\eqref{s-eq}-\eqref{ssss-eq} somewhat later.
\begin{align*}
\s^i_jx&\de c_i(d_{ij}\cdot x)\quad\text{for\ }i\ne j, \\
\s^{01}_{12}x&\de \s^0_1\s^1_2x,    \\
\s^{01}_{10}x&\de {}_2s(0,1)c_2x = \s^2_0\s^0_1\s^1_2c_2x,   \\
\s^{12}_{01}x&\de  \s^2_1\s^1_0x   .
\end{align*}
Next we write up the terms we use in expressing that it is not the
case that the binary relations at places $01$ and $12$ of $x$
coincide:
\begin{align*}
\gamma(x)\de c_2x\oplus s^{12}_{01}c_0x .
\end{align*}
Finally, the terms for expressing that $c_2x$ is not a strict linear
order:
\begin{align*}
\iota(x)&\de c_2x\cdot d_{01}\quad&&\mbox{not irreflexive,}\\
\sigma(x)&\de c_2x\cdot s^{01}_{10}c_2x\quad&&\mbox{not antisymmetric,}\\
\tau(x)&\de c_2x\cdot s^{01}_{12}c_2x-s^1_2c_2x\quad&&\mbox{not transitive,}\\
\lambda(x)&\de c_1c_2x\cdot c_0c_2x\cdot -c_2x\cdot -s^{01}_{10}c_2x\quad&&\mbox{not linear,}\\
\o(x)&\de \iota(x)+\sigma(x)+\tau(x)+\lambda(x) .&
\end{align*}
Let
\begin{equation}\label{cutt-eq}
e(x)\ \ \de \ \ x\le c_{(3)}(\beta(x)+\gamma(x)+\o(x)) ,
\end{equation}
where $c_{(3)}y\de c_0c_1c_2y$.

We now turn to stating precisely what the equation $e(x)$ expresses
in a set-algebra about an element $x$.
We will use the following notation extensively. Let $U$ be a set,
$s\in {}^\a U$, let $n\in\omega$, let $H\in{}^n\a$ be
repetition-free, and let $q\in{}^nU$. Then $s(H/q)$ denotes the
sequence we get from $s$ by changing $s(H_k)$ to $q_k$,
simultaneously, for all $k<n$. We will write finite sequences
$\langle i_0,i_1,...,i_n\rangle$ in the simplified form
$i_0i_1...i_n$ when this is not likely to lead to confusion. E.g.,
$s(12/uv)=s(1/u)(2/v)$. Further, if $x\subseteq {}^\a U$, then
$x[s,H]$ denotes the $n$-place relation defined as
\begin{align*}
x[s,H]&\de \{ q\in{}^nU : s(H/q)\in x\}.
\end{align*}
For example, $x[s,01]=\{ uv : s(01/uv)\in x\}$, and $x[s,0]=\{ u :
s(0/u)\in x\}$.

\begin{lem}\label{express-lem}
Let $\Cc$ be any $\a$-dimensional set algebra with unit a disjoint
union of weak spaces.
 Then the equation $e$ as defined in \eqref{cutt-eq} is
true in $\Cc$ at $x\in C$ iff for all $s\in x$ it is true that
either $x[s,012]$ is not a sensitive cut of $Z\de (c_0x\cdot
c_2x)[s,012]$, or $<_x\de c_2x[s,01]$ is not a strict linear order
on $W\de c_1c_2x[s,0]$, or else $Z\ne \{ uvw\in {}^3W : u<_x v<_x
w\}$.
\end{lem}

\noindent {\bf Proof.} Let $z\de c_0x\cdot c_2x$ and $s\in x$.
\begin{equation}\label{bigcut-eq}
 s\notin c_{(3)}\beta(x) \ \ \mbox{ iff }\ \ x[s,012] \mbox{ is a
sensitive cut of }z[s,012].
\end{equation}
Indeed,
\begin{alignat*}{2}
&s\in-c_{(3)}(c_0x\oplus c_0(z-x)) &\quad& \mbox{iff}\\
&s(012/uvw)\notin c_0x\oplus c_0(z-x)\mbox{ for all $u,v,w$}  &\quad& \mbox{iff}\\
&s(012/uvw)\in (c_0x\cdot c_0(z-x)+(-c_0x\cdot-c_0(z-x))\mbox{ for all $u,v,w$\quad}  &\quad& \mbox{iff}\\
&(c_0x)[s,012]=(c_0(z-x))[s,012]\mbox{}  &\quad& \mbox{iff}\\
&c_0(x[s,012])=c_0(z[s,012]-x[s,012]) .\mbox{}  &\quad& \mbox{}
\end{alignat*}
In the last step we used $(c_0x)[s,012]=c_0(x[s,012])$ and
$(z-x)[s,012] = z[s,012] - x[s,012]$ which statements are easy to
verify by using the definitions.
All the above hold also for $c_1$ and $c_2$ in place of $c_0$, so we
get \eqref{bigcut-eq}.

For dealing with the rest of the terms, we make some preparations.
We will check the following: Assume $i\ne j$.
\begin{alignat}{3}
\label{s-eq}
&s\in \s^i_jx       &\qquad& \mbox{ iff } &\qquad& s(i/s_j)\in x, \\
\label{ss-eq}
&s\in \s^{12}_{01}x &\qquad& \mbox{ iff } &\qquad& s(12/s_0s_1)\in x,    \\
\label{sss-eq}
&s\in \s^{01}_{10}c_2x &\qquad& \mbox{ iff } &\qquad& s(01/s_1s_0)\in c_2x,    \\
\label{ssss-eq}
&s\in \s^{01}_{12}x &\qquad& \mbox{ iff } &\qquad&
s(01/s_1s_2)\in x.
\end{alignat}
Indeed, \eqref{s-eq} is true because
\begin{alignat*}{2}
&s\in \s^i_jx = c_i(d_{ij}\cdot x) &\quad& \mbox{iff}\\
&s(i/u)\in d_{ij}\cdot x\quad\mbox{\ for some $u$ } &\quad& \mbox{iff}\\
&u=s_j\mbox{\ and\ }s(i/u)\in x\quad\mbox{\ for some $u$ } &\quad& \mbox{iff}\\
&s(i/s_j)\in x .\quad\mbox{} &\quad& \mbox{}
\end{alignat*}
\eqref{ss-eq} is true because
\begin{alignat*}{2}
&s\in \s^{12}_{01}x = \s^2_1\s^1_0x &\quad& \mbox{iff, by \eqref{s-eq}}\\
&s(2/s_1)\in \s^1_0x\quad\mbox{} &\quad& \mbox{iff, by \eqref{s-eq}}\\
&(s(2/s_1))(1/s_0)\in x\quad\mbox{} &\quad& \mbox{iff}\\
&s(12/s_0s_1)\in x .\quad\mbox{} &\quad& \mbox{}
\end{alignat*}
Checking \eqref{sss-eq}:
\begin{alignat*}{2}
&s\in \s^{01}_{10}c_2x = \s^2_0\s^0_1\s^1_2c_2x &\quad& \mbox{iff, by \eqref{s-eq}}\\
&s(2/s_0)\in \s^0_1\s^1_2c_2x\quad\mbox{} &\quad& \mbox{iff, by \eqref{s-eq}}\\
&(s(2/s_0))(0/s_1)\in \s^1_2c_2x\quad\mbox{} &\quad& \mbox{iff, by \eqref{s-eq}}\\
&((s(2/s_0))(0/s_1))(1/s_0)\in c_2x\quad\mbox{} &\quad& \mbox{iff}\\
&s(012/s_1s_0s_0)\in c_2x\quad\mbox{} &\quad& \mbox{iff}\\
&s(01/s_1s_0)\in c_2x .\quad\mbox{} &\quad& \mbox{}
\end{alignat*}
Checking \eqref{ssss-eq}:
\begin{alignat*}{2}
&s\in \s^{01}_{12}x = \s^0_1\s^1_2x &\quad& \mbox{iff, by \eqref{s-eq}}\\
&s(0/s_1)\in \s^1_2x\quad\mbox{} &\quad& \mbox{iff, by \eqref{s-eq}}\\
&(s(0/s_1))(1/s_2)\in x\quad\mbox{} &\quad& \mbox{iff}\\
&s(01/s_1s_2)\in x .\quad\mbox{} &\quad& \mbox{}\\
\end{alignat*}
We are ready to continue with the terms occurring in $e$.%
\begin{equation}\label{closure-eq}
 s\notin c_{(3)}\gamma(x) \ \ \mbox{ iff }\ \ c_2x[s,01]=c_0x[s,12].
\end{equation}
Indeed,
\begin{alignat*}{2}
&s\in-c_{(3)}\gamma(x)=-c_{(3)}(c_2x\oplus \s^{12}_{01}c_0x) &\quad& \mbox{iff}\\
&s(012/uvw)\notin c_2x\oplus \s^{12}_{01}c_0x\mbox{\ \ for all $u,v,w$} &\quad& \mbox{iff}\\
&s(012/uvw)\in (c_2x\cdot \s^{12}_{01}c_0x+(-c_2x\cdot-\s^{12}_{01}c_0x)\mbox{\ \  for all $u,v,w$\quad} &\quad& \mbox{iff}\\
&s(012/uvw)\in c_2x\mbox{ iff }s(012/uvw)\in\s^{12}_{01}c_0x\mbox{\ \  for all $u,v,w$} &\quad& \mbox{iff, by \eqref{ss-eq}} \\
&uv\in c_2x[s,01]\mbox{ iff }s(012/wuv)\in c_0x\mbox{\ \  for all $u,v,w$} &\quad& \mbox{iff} \\
&uv\in c_2x[s,01]\mbox{ iff }uv\in c_0x[s,12],\mbox{\ \  for all $u,v$}&\quad& \mbox{iff} \\
&c_2x[s,01]=c_0x[s,12] .&\quad& \mbox{}
\end{alignat*}
By this, \eqref{closure-eq} has been proved. We note that
$c_2x[s,01]=c_0x[s,12]$ implies that
$\Rg(c_2x[s,01])=
\Dom(c_2x[s,01])$.  Indeed,
\begin{alignat*}{3}
&u\in \Rg(c_2x[s,01])&\quad& &\qquad\quad& \mbox{iff}\\
&vu\in c_2x[s,01]&\quad&\text{ for some }v &\qquad\quad& \mbox{iff}\\
&s(012/vuw)\in x&\quad&\mbox{ for some }v,w &\qquad& \mbox{iff}\\
&s(012/vuw)\in c_0x&\quad&\mbox{ for some }v,w &\qquad& \mbox{iff}\\
&uw\in c_0x[s,12]&\quad&\mbox{ for some }w &\qquad& \mbox{iff}\\
&u\in\Dom(c_0x[s,12])&\quad& &\qquad& \mbox{iff, by $c_0x[s,12]=c_2x[s,01]$}\\
&u\in \Dom(c_2x[s,01]) .& &\qquad\quad& \mbox{}\\
\end{alignat*}

We say that a binary relation $R$ is \emph{linear on $W$} iff $W$ is
both the domain and range of $R$ and $\langle u,v\rangle\in R$ or
$\langle v,u\rangle\in R$ for all $u,v\in W$. We have seen that
$s\in-c_{(3)}\gamma(x)$ implies that the domain of $c_2x[s,01]$
coincides with its range. Assume for $\eqref{order-eq}$ below that
the domain and range of $c_2x[s,01]$ coincide. We can do this
because we will use \eqref{order-eq} only when $s\notin
c_{(3)}\gamma(x)$.

\begin{equation}\label{order-eq}
s\notin c_{(3)}\o(x)   \ \mbox{ iff }\  c_2x[s,01] \mbox{ is a
strict linear order on }c_1c_2x[s,0].
\end{equation}
Indeed,
\begin{alignat*}{2}
&s\in-c_{(3)}\iota(x)=-c_{(3)}(c_2x\cdot d_{01}) &\quad& \mbox{iff}\\
&s(012/uvw)\notin (c_2x\cdot d_{01})\mbox{ for all $u,v,w$} &\quad& \mbox{iff}\\
&s(012/uvw)\in c_2x\mbox{ implies }u\ne v\mbox{ for all $u,v,w$\quad} &\quad& \mbox{iff}\\
&u\ne v \mbox{ for all $uvw\in c_2x[s,012]$} &\quad& \mbox{iff} \\
&u\ne v \mbox{ for all $uv\in c_2x[s,01]$} &\quad& \mbox{iff} \\
&c_2x[s,01]\text{ is irreflexive.}
\end{alignat*}
\begin{alignat*}{2}
&s\in-c_{(3)}\sigma(x)=-c_{(3)}(c_2x\cdot \s^{01}_{10}c_2x) &\quad& \mbox{iff}\\
&s(012/uvw)\notin (c_2x\cdot \s^{01}_{10}c_2x)\mbox{\ \ for all $u,v,w$} &\quad& \mbox{iff}\\
&s(012/uvw)\in c_2x\mbox{ implies }s(012/uvw)\notin \s^{01}_{10}c_2x\quad &\quad& \mbox{iff, by \eqref{sss-eq}}\\
&uv\in c_2x[s,01]\mbox{ implies }vu\notin c_2x[s,01] &\quad& \mbox{iff}\\
&c_2x[s,01]\mbox{ is antisymmetric.}
\end{alignat*}
\begin{alignat*}{2}
&s\in-c_{(3)}\tau(x)=-c_{(3)}(c_2x\cdot \s^{01}_{12}c_2x-\s^1_2c_2x) &\quad& \mbox{iff}\\
&s(012/uvw)\notin (c_2x\cdot s^{01}_{12}c_2x-s^1_2c_2x)\mbox{\ \  for all $u,v,w$} &\quad& \mbox{iff}\\
&s(012/uvw)\in c_2x\cdot \s^{01}_{12}c_2x\mbox{ implies }s(012/uvw)\in \s^1_2c_2x\quad &\quad& \mbox{iff, by \eqref{ssss-eq},\eqref{s-eq}}\\
&uv,vw\in c_2x[s,01]\mbox{ implies }uw\in c_2x[s,01] &\quad& \mbox{iff}\\
&c_2x[s,01]\mbox{ is transitive.}&\quad& \mbox{}
\end{alignat*}
\begin{alignat*}{2}
&s\in-c_{(3)}\lambda(x)=-c_{(3)}(c_1c_2x\cdot c_0c_2x\cdot -c_2x\cdot -\s^{01}_{10}c_2x) &\quad& \mbox{iff}\\
&s(012/uvw)\notin (c_1c_2x\cdot c_0c_2x\cdot -c_2x\cdot -\s^{01}_{10}c_2x)\mbox{ for all $u,v,w$} &\quad& \mbox{iff}\\
&s(012/uvw)\in c_1c_2x\cdot c_0c_2x\mbox{ $\Rightarrow$ }s(012/uvw)\in c_2x+\s^{01}_{10}c_2x &\quad& \mbox{iff, by \eqref{sss-eq}}\\
&u\in \Dom(c_2x[s,01]), v\in \Rg(c_2x[s,01])\mbox{ $\Rightarrow$ }uv\in c_2x[s,01]\mbox{ or }vu\in c_2x[s,01]&\quad& \mbox{iff}\\
&c_2x[s,01]\text{ is linear on its field if
$\Dom(c_2x[s,01])=\Rg(c_2x[s,01])$ .}&\quad& \mbox{}\\
\end{alignat*}

By the above, \eqref{order-eq} has been proved. Let $x\in C$. By
\eqref{cutt-eq}, then $\Cc\models e(x)$ iff for all $s\in x$ we have
$s\in c_{(3)}(\beta(x)+\gamma(x)+\o(x))$. Now, $s\in
c_{(3)}(\beta(x)+\gamma(x)+\o(x))$ iff $s\in c_{(3)}\beta(x)$ or
$s\in c_{(3)}\gamma(x)$ or we have $s\in c_{(3)}\o(x)$ when $s\notin
c_{(3)}\gamma(x)$. By \eqref{order-eq}, \eqref{closure-eq},
\eqref{bigcut-eq}, then $s\in c_{(3)}(\beta(x)+\gamma(x)+\o(x))$ iff
it is not the case that $<_x=c_2x[s,01]$ is a strict linear order,
$Z=c_0x\cdot c_2x[s,012]=\{ \langle u,v,w\rangle : u<_x v<_x w \}$
and $x[s,012]$ is a sensitive cut of $Z[s,012]$. \QED

\subsection{Checking the equations in the witness algebra}\label{check-sec}

\begin{lem}\label{e-fails} The equation $e$ fails in $\Aa$.
\end{lem}

\noindent {\bf Proof.} We show that $e$ fails in $\Aa$ at $g$. Let
$s\in g$ be arbitrary such that $s$ agrees with $p$ on all indices
$i>2$. There is such a sequence. Then, $<_g=c_2g[s,01]=c_2T[s,01]=\{
\langle u,v\rangle : u,v\in V_0, u<v\}$ is a strict linear order on
$c_1c_2g[s,0]=c_1c_2T[s,0]=V_0$. Also, $Z=(c_0g\cdot
c_2g)[s,012]=T[s,012]=\{ \langle u,v,w\rangle\in {}^3V_0 : u<v<w\}$.
Finally,  $g[s,012]=\{\langle u, (u+v)/2, v\rangle : u<v, u,v\in
V_0\}$ is a sensitive cut of $Z$. Then Lemma~\ref{express-lem}
implies that $e$ fails in $\Aa$ at $g$. \QED

Let $i,j,k\in\a-\{0,1,2\}$ be distinct and let $e_{ijk}$ denote the
equation we get from $e$ by replacing the indices $0,1,2$ everywhere
with $i,j,k$ respectively. We are going to show that $e_{ijk}$ holds
in $\Aa$.
Lemma~\ref{express-lem} is true with systematically replacing the
indices $0,1,2$ by $i,j,k$. Thus, to show that $e_{ijk}$ holds in
$\Aa$, we have to show that $x[s,ijk]$ is not a sensitive cut of a
ternary relation built up from a linear order $<_x$ in the way $T$
is built up from $<$, for all $x\in A$ and $s\in x$.
In proving this, we will use the following lemma, which says that
certain permutations of $U$ leave the relations $x[s,ijk]$, that
determine the validity of $e_{ijk}$, fixed. We agree on some
terminology first.

Let $\pi$ be a permutation of $U$, and let $n$ be an ordinal. Then
$\pi(s)\de\langle \pi(s_i) : i<n\rangle$ if $s\in{}^nU$, and
$\pi(R)\de\{\pi(s) : s\in R\}$ for $R\subseteq{}^nU$. We say that
$\pi$ leaves $R$ fixed iff $\pi(R)=R$. In the present section we
shall often use a certain property of permutations of $U$, so we
give it a temporary name.

\begin{defi}\label{good-def}
We say that $\pi$ is \emph{good} iff it satisfies {\rm(i)-(iii)}
below.
\begin{description}
\item{\rm (i)}
$\pi$ leaves $<$ fixed, i.e., $u<v$ iff $\pi(u)<\pi(v)$ for all
$u,v\in V_0$,
\item{\rm(ii)}
$\pi$ leaves all the $V_m$s fixed, i.e., $\pi(V_m)=V_m$ for all
$m\in\omega$,
\item{\rm(iii)}
$\pi$ is the identity on all but a finite number of $V_m$s, i.e.,
$\{ m\in\a : \text{$\forall u\in V_m(\pi(u)=u)$}\}$ is a co-finite
subset of $\a$.
\end{description}
\end{defi}

\begin{lem}\label{aut-lem}
For any $x\in A$ and $s\in V$ there is a finite $S\subseteq U$ such
that any good permutation of $U$ which is identity on $S$ leaves
$x[s,ijk]$ fixed.
\end{lem}

\noindent{\bf Proof.} Let $x\in A$. Then $x=y+h$ for some $y$
generated by $T$ and for some $h\in\{ 0,g,T-g\}$, by
Lemma~\ref{gen-lem}. Assume that $y=\xi(T)$ for a term $\xi$. Let
$\Delta\subseteq\a$ be finite such that it contains all the indices
occurring in $\xi$ as well as $0,1,2$. We show that for all $s,s'$
\begin{equation}\label{iff-eq}
s\in T\leftrightarrow s'\in T,\ \
s\upharpoonright\Delta=s'\upharpoonright\Delta\quad\ \text{ imply
}\quad \ s\in y\leftrightarrow s'\in y.
\end{equation}
We prove \eqref{iff-eq} by induction on elements $z$ generated from
$T$ by the use of indices from $\Delta$. Clearly, \eqref{iff-eq}
holds for $T$ and $d_{mn}$ for $m,n\in\Delta$. Assume that
\eqref{iff-eq} holds for $z,z'$. Then clearly it holds for $-z$ and
$z\cdot z'$. Let $m\in\Delta$ and assume that $s,s'$ satisfy the
conditions. Now,
\begin{alignat*}{2}
&s\in c_mz&\qquad&\text{iff, by the definition of $c_m$}\\
&s(m/u)\in z\text{ for some $u$}&\qquad&\text{iff, by the induction hyp., see details below}\\
&s'(m/u)\in z\text{ for the same $u$}&\qquad&\text{iff, by the definition of $c_m$}\\
&s'\in c_mz .&\qquad&\text{}
\end{alignat*}
Above, in the step from the second to third line we used that
$s'(m/u)$ agrees with $s(m/u)$ on $\Delta$ and $s(m/u)\in T$ iff
$s'(m/u)\in T$ (by $s\in T$ iff $s'\in T$, the definition of $T$,
and $0,1,2\in\Delta$). By this, \eqref{iff-eq} has been proved.
\smallskip

Let now $s\in V$ and
\[S\de \{ s_m : m\in\Delta\}\cup V_i\cup V_j\cup V_k. \]
Then $S\subseteq U$ is finite. Let $\pi$ be a good permutation of
$U$ which is the identity on $S$, we want to show that $\pi$ leaves
$x[s,ijk]$ fixed. Recall that $x=y+h$ where $h\in\{0,g,T-g\}$. Then
$x[s,ijk]=y[s,ijk]+h[s,ijk]$. Thus
$\pi(x[s,ijk])=\pi(y[s,ijk])+\pi(h[s,ijk])$, so it is enough to show
that $\pi$ leaves both $y[s,ijk]$ and $h[s,ijk]$ fixed. We begin
with the second.
Indeed, $h[s,ijk]\subseteq V_i\times V_j\times V_k$ when
$h\in\{0,g,T-g\}$, thus $\pi(uvw)=uvw$ for all $uvw\in h[s,ijk]$ by
$V_i\cup V_j\cup V_k\subseteq S$ and $\pi$ being the identity on
$S$.

We turn to $y[s,ijk]$. First we note that $\pi$ being good implies
that $\pi(V)=V$ by Def.\ref{good-def}(iii), and then $\pi(T)=T$ by
Def.\ref{good-def}(i),(ii). Since $y$ is generated by $T$ we then
have (by, e.g., \cite[3.1.36]{HMTII})
\begin{equation}\label{pi-eq}
\pi(y)=y .
\end{equation}
We want to show that
\begin{equation}\label{y-eq}
uvw\in y[s,ijk]\quad\text{ iff }\quad\pi(uvw)\in y[s,ijk].
\end{equation}
Indeed,
\begin{alignat*}{2}
&uvw\in y[s,ijk]&\qquad&\text{iff, by the definition of $y[s,ijk]$}\\
&s(ijk/uvw)\in y&\qquad&\text{iff, by \eqref{pi-eq}} \\
&\pi(s(ijk/uvw))\in y&\qquad&\text{iff, by \eqref{iff-eq} and see below}\\
&s(ijk/\pi(uvw))\in y&\qquad&\text{iff, by the definition of $y[s,ijk]$}\\
&\pi(uvw)\in y[s,ijk] .&\qquad&\text{}
\end{alignat*}
In the argument from the third to fourth line we used that
$\pi(s(ijk/uvw))$ and $s(ijk/\pi(uvw))$ agree on $i,j,k$ by their
definitions, they agree on $\Delta-\{ i,j,k\}$ by $\pi$ being the
identity on $S\supseteq\{ s_m : m\in\Delta\}$; further, one of them
is in $T$ iff the other is:
\begin{alignat*}{2}
&\pi(s(ijk/uvw))\in T&\qquad&\text{iff, by $\pi(T)=T$}\\
&s(ijk/uvw)\in T&\qquad&\text{iff, by Def.\ref{good-def}(ii)} \\
&s(ijk/\pi(uvw))\in T .&\qquad&\text{}
\end{alignat*}
This proves \eqref{y-eq}, and Lemma~\ref{aut-lem} has been
proved.\QED

\begin{lem}\label{check-lem}
The equation $e_{ijk}$ is valid in $\Aa$ when $i,j,k\in\a-\{0,1,2\}$
are distinct.
\end{lem}

\noindent {\bf Proof.} For checking validity of $e_{ijk}$, we will
use Lemma~\ref{express-lem} (with $0,1,2$ systematically replaced by
$i,j,k$). Let $x\in A$, $s\in x$ and $R\de x[s,ijk]$, $W\de
c_jc_kx[s,i]$, $Z\de (c_ix\cdot c_kx)[s,ijk]$. Assume that $<_x\de
c_kx[s,ij]$ is a strict linear order on $W$ and $Z=\{\langle
u,v,w\rangle\in{}^3W : u<_x v<_x w\}$. We have to show that $R$ is
not a sensitive cut of $Z$. Let $S\subseteq U$ be such that
\begin{equation}\label{S-eq}
\text{$\pi(R)=R$ for all good permutations $\pi$ of $U$ that are
identity on $S$.} \end{equation}
There is such an $S$ by Lemma~\ref{aut-lem}. We note that
\begin{equation}\label{o-eq}
\text{$\pi(R)=R$\quad implies\quad $\pi(<_x)=<_x$ and $\pi(W)=W$ .}
\end{equation}
This is true because $<_x=c_2R$ by their definitions:
$<_x=c_kx[s,ij]=\{uv : s(ij/uv)\in c_kx\} = c_2\{ uvw :
s(ijk/uvw)\in x\} = c_2R$. Below,  we will use (*) from the proof of
Lemma~\ref{gen-lem} several times.

We turn to showing that $R$ is not a sensitive cut of $Z$. By $s\in
x$ we have that $s_is_j\in c_kx[s,ij]$ so $<_x$ is nonzero. Since by
our assumption $<_x$ is a strict linear order on $W$, it does not
have a maximal element (by $W=\Dom{(<_x)}=\Rg(<_x)$), so $W$ is
infinite by $W\ne\emptyset$. Assume $m\ge 3$ is such that $S$ is
disjoint from $V_m$. We show that $V_m$ is disjoint from $W$. Assume
$W\cap V_m\ne\emptyset$.
Let $\pi$ be the permutation of $U$ that interchanges the elements
of $V_m$ and it leaves all the other elements of $U$ fixed. (Recall
that the $V_m$'s for $m\ge 3$ have two elements.) Then $\pi$ is good
and it is identity on $S$, so it leaves $W$ as well as $<_x$ fixed,
by \eqref{S-eq},\eqref{o-eq}. This implies that $V_m\subseteq W$ by
$V_m\cap W\ne\emptyset$, so by $<_x$ being linear on $W$, we have
$a<_x b$ for some distinct $a,b\in V_m$. By $\pi$ leaving $<_x$
fixed, then we have $b<_x a$ (by $\pi(a)=b$, $\pi(b)=a$). This
contradicts $<_x$ being antisymmetric.

Thus $W$ intersects only finitely many of the $V_m$s. Then $W\cap
V_0$ is infinite because all the $V_m$s disjoint from $V_0$ are
finite. Let
\[ K\de W\cap V_0\cap S\qquad\text{and}\qquad W'\de (W\cap V_0)-S .\]
Thus $K$ is finite and $W'$ is infinite. Therefore, there are
distinct $u,v\in W'$ such that no element of $K$ lies in between
$u,v$ according to $<$. (Indeed, let $K=\{ k_1,\dots , k_n\}$ with
$k_1<\dots < k_n$. Then there are at least two elements of $W'$ that
lie in the same interval determined by the $k_m$s, and they will
do.) We may assume $u<v$.
Since $<_x$ is linear on $W$ and $u,v\in W$, we have either $u<_x v$
or $v<_x u$. We assume $u<_x v$, the case $v<_x u$ will be
completely analogous, see \eqref{betw2-eq}. So assume
\[ [u,v]\cap K=\emptyset, \quad u<v ,\quad\text{\ and\ }\quad u<_x v \]
and we are going to show that for all $w\in U$ we have
\begin{align}\label{betw-eq}
u<_x w <_x v\qquad\text{iff}\qquad  u < w < v .
\end{align}
Indeed, to prove \eqref{betw-eq}, let first $u < w < v $ be
arbitrary, we want to show $u<_x w <_x v$. Let $u'<u$ be such that
there is no element of $K$ between $u'$ and $u$, there is such an
$u'$ because $<$ is dense and $K$ is finite. Then there is no
element of $K$ between $u'$ and $v$. Take a $\pi$ as in (*) for
$u'<u<v$ and $u'<w<v$ and extend it to $U$ by being the identity on
$U-V_0$. Then this $\pi$ is identity on $S$ and it is good. So it
leaves $<_x$ fixed, by \eqref{S-eq},\eqref{o-eq}. By $u<_x v$ then
we have $w=\pi(u)<_x \pi(v)=v$. By a similar argument we get $u <_x
w$. (Indeed, choose $v'>v$ such that there is no element of $K$
between $v$ and $v'$ and apply (*) with $u<v<v'$ and $u<w<v'$.)
We have seen $u <_x w <_x v$. To prove the other direction, assume
that $w\in U$ and it is not the case that $u < w < v$. Thus either
$w<u<v$ or $u<v<w$.
In either cases, there is a good permutation $\pi$ of $U$ which is
identity on $S$, leaves $w$ fixed and takes $u$ to $v$, and there is
also a good permutation $\pi$ of $U$ which is identity on $S$,
leaves $w$ fixed and takes $v$ to $u$. (Indeed, take $u'<u$ and
$v<v'$ such that no element of $K\cup\{ w\}$ lies between $u'$ and
$v'$, and then apply (*).)
Hence $w<_x u$ iff $w<_x v$  and   $u<_x w$ iff $v<_x w$. Hence, it
is not the case that $u<_x w <_x v$, as it was desired. The equation
\eqref{betw-eq} has been proved.\smallskip

Assume now the other case, i.e., that
\[ [u,v]\cap K=\emptyset, \quad u<v ,\quad\text{\ and\ }\quad v<_x u. \]
We are going to show that for all $w\in U$ we have
\begin{align}\label{betw2-eq}
v<_x w <_x u\qquad\text{iff}\qquad  u < w < v .
\end{align}
To prove \eqref{betw2-eq}, let first $u < w < v $ be arbitrary, we
want to show $u<_x w <_x v$. Let $u'<u$ be such that there is no
element of $K$ between $u'$ and $u$. Take a $\pi$ as in (*) for
$u'<u<v$ and $u'<w<v$ and extend it to $U$ by being the identity on
$U-V_0$. Then this $\pi$ is identity on $S$ and it is good. So it
leaves $<_x$ fixed, by \eqref{S-eq},\eqref{o-eq}. By $v<_x u$ then
we have $v=\pi(v)<_x \pi(u)=w$. We get $w <_x u$ by choosing $v'>v$
such that there is no element of $K$ between $v$ and $v'$ and
applying (*) with $u<v<v'$ and $u<w<v'$. We have seen $v <_x w <_x
u$. The proof of the other direction is the same as in the proof for
\eqref{betw-eq}. The equation \eqref{betw2-eq} has been
proved.\smallskip

We are ready to prove that $R$ is not a sensitive cut of $Z$. Assume
that $Z\subseteq c_1R$, we will show that $Z\not\subseteq c_1(Z-R)$.
By \eqref{betw-eq},\eqref{betw2-eq} and $<$ being dense we have that
$uwv\in Z$ for some $w$ (and the $u,v$ chosen as before), so
$uw'v\in R$ for some $w'$ by $Z\subseteq c_1R$. Let $u<w''<v$ be
arbitrary and take a good permutation $\pi$ of $U$ that takes $w'$
to $w''$ and leaves everything outside the open interval $(u,v)$
fixed. There is such a $\pi$ by (*). This $\pi$ leaves $S$ fixed
since no element of $S$ lies between $u$ and $v$ (according to $<$).
Then it leaves $R$ fixed by \eqref{S-eq}. So $uw''v\in R$ by
$uw'v\in R$ and $\pi(uv'v)=uv''v$. By
\eqref{betw-eq},\eqref{betw2-eq} this means that $R(uw'v)$ for all
$w'$ such that $uw'v\in Z$. Hence $uwv\in Z$ is such that $uwv\notin
c_1(Z-R)$, and we are done with proving Lemma~\ref{check-lem}.\QED

We now round up the proof of Theorem~\ref{manyvar-thm}. Let
$I\subseteq\a\times\a\times\a$ be such that $|I|=|\a|$ and for all
distinct $ijk, lmn\in I$ we have $\{ i,j,k\}\cap\{
l,m,n\}=\emptyset$, $|\{ i,j,k\}|=3$ and $012\in I$. There is such
an $I$ since $\a$ is infinite. Let
\[ E \de \{ e_{ijk} : ijk\in I \} .\]
Then $|E|=|\a|$. For all $ijk\in I$ let $\Aa_{ijk}$ denote the
algebra we get from $\Aa$ by interchanging (renaming) the operations
$c_m, d_{mn}$ for $m,n\in \{0,1,2\}$ with those for $m,n\in
\{i,j,k\}$, respectively.
This algebra is denoted by $\Rd^{\rho}\Aa$ where $\rho$ is the
permutation of $\a$ which interchanges $i,j,k$ with $0,1,2$
respectively, see \cite[Def.2.6.1]{HMTI}. By
Lemmas~\ref{e-fails},\,\ref{check-lem} then we have
\begin{equation}\label{rd-eq}
\Aa_{ijk}\not\models e_{ijk}\quad\mbox{ while }\quad\Aa_{ijk}\models
e_{klm}\ \mbox{ for all $klm\in I, klm\ne ijk$}.
\end{equation}
Also, $\Aa_{ijk}\in\Iso\Ws_\a$ by \cite[3.1.119]{HMTII} and, so
$\Aa_{ijk}\in\iCs_\a$ by \cite[3.1.102]{HMTII} and because the fact
that $\Aa$ has an infinite base is reflected on its equational
theory.
(We note that $\Aa_{ijk}$ is isomorphic to the algebra we get from
$\Aa$ by replacing in its construction $0,1,2$ with $i,j,k$
systematically, and changing nothing else. An isomorphism showing
this takes $x\in\Aa_{ijk}$ to $\{\rho(s) : s\in x\}$.)

For $G\subseteq E$ define
\[ \V_G\de\{ \Bb\in\iCs_\a : \Bb\models G\} .\]
Then $\V_G$ is a subvariety of $\iCs_\a$. Assume $G,H\subseteq E$
are distinct. Then there is $ijk\in I$ distinguishing them, we may
assume $ijk\in G$ and $ijk\notin H$. By \eqref{rd-eq} we have that
$\Aa_{ijk}\notin\V_G$ but $\Aa_{ijk}\in\V_H$, so $\V_G\ne\V_H$.
Therefore, there are $2^{|\a|}$ distinct subvarieties of $\iCs_\a$.
The same is true for $\RCA_\a$ by $\iCs_\a\subseteq\RCA_\a$.
Theorem~\ref{manyvar-thm} has been proved.\bigskip

We close this section with discussing some properties necessary for
our construction $\Aa$ to work.

\begin{rem}\label{constr-rem}
\rm{ (i) It is necessary that the base of $\Aa$ (i.e., the set $U$)
be at least of cardinality $|\a|$. This is true because algebras of
smaller base are diagonal (roughly: each of their elements
intersects many diagonal elements, for precise definition see
\cite[p.416]{HMTI}), and we will prove that all diagonal algebras
are symmetric, see Theorem~\ref{sym-thm}. Clearly, $\Aa$ has to be
non-symmetric to play its role in the proof.

(ii) In the equation $e$ it was necessary to code a property that
can occur on an infinite set only, this is the role of using the
ordering on rational numbers $V_0=V_1=V_2$ in the definition of the
generator element $g$. In more detail: an equation  $e(x)$ using
indices from $\{ i,j,k\}$ can talk in a set algebra about the
ternary relation $x[s,ijk]$ only. However, all ternary relations on
a finite set occur as $x[s,ijk]$ in a set algebra when the base set
is infinite. (This is the main idea used in \cite{NAPAL87}.) Since
we want $e_{012}(x)$ to hold and $e_{ijk}(x)$ to fail in our witness
algebra, $e(x)$ has to code a property of $x[s,ijk]$ which can be
realized only on infinite sets.

(iii) The equation $e$ fails in $\Aa$ at $g$, but $e$ is true in
$\Aa$ for all elements that are closed to at least one
cylindrification $c_i$. Indeed, we can see that $e$ holds for
$c_ia\in A$ as follows. Lemma~\ref{gen-lem} and \eqref{cut-eq}
together with $T\in B$ imply that $c_ia\in B$. The proof of
Lemma~\ref{aut-lem} works for $x\in B$ and $ijk=012$, and then the
proof of Lemma~\ref{check-lem} works to show that $e$ holds for
$c_ia$. Thus, in $\Aa$ cylindrification-closed elements satisfy more
equations than all the elements. This behavior of our witness
algebra $\Aa$ is necessary, because each algebra in which no such
behavior occurs is symmetric (see Theorem~\ref{ret-thm}(i) in
subsection~\ref{ind-ssec}). }
\end{rem}

\section{Subvarieties of $\CA_\a$ containing $\RCA_\a$}\label{ca-sec}

This section contains an unpublished theorem from \cite{NVarprep86}.
The proof is analogous to the proof of Theorem~\ref{manyvar-thm}.

\begin{thm}\label{ca-thm} There are $2^{|\a|}$ distinct subvarieties of
$\CA_\a$ all containing $\RCA_\a$.
\end{thm}

\noindent {\bf Proof.} We are going to exhibit an equation $e$ valid
in $\RCA_\a$ and an algebra $\Aa\in\CA_\a$ such that $\Aa\not\models
e$ while $\Aa\models e_{i1}$ for appropriate versions $e_{i1}$ of
$e$. This $e$ is Henkin's equation $e_{ij}(x,y)$ with $ij$ taken as
$01$:
\begin{equation}
e_{ij}(x,y)\quad\de\quad c_j(x\cdot y\cdot c_i(x-y))\le
c_i(c_jx-d_{ij}),
\end{equation}
see \cite[3.2.65]{HMTII}. For a simplified version of this equation
see \cite[chap.3.5]{Venema1}, and for a drawing see
\cite[p.551]{Nsurv}. Henkin's equation expresses that if the domains
of $R$ and $S$ coincide and this common domain is a singleton, then
$R$ and $S$ are disjoint iff their ranges are disjoint. Now,
$\RCA_\a\models e_{ij}(x,y)$, by e.g. \cite[3.2.65]{HMTII}.

We now turn to constructing our ``witness" algebra $\Aa$. It is
obtained from a representable algebra $\Bb$ in which we split an
atom whose domain is a singleton into two parts both having the same
domain and range as the original atom. Henkin's equation then will
fail for the split elements. In some sense this will be a
``nonrepresentable counterpart" of the construction we used in the
proof of Theorem~\ref{manyvar-thm}.

Let $\langle V_i : i\in\a\rangle$ be a system of disjoint sets such
that $V_0$ is a singleton, and $V_i$ for $i\ge 1$ have more than one
elements. Let $U$ be the union of these sets, let $V\de{}^\a U$ and
let $g$ be the direct product of the $V_i$,  i.e.,
\[ g \de \prod\langle V_i : i\in\a\rangle \de \{ s\in V : s_i\in
V_i\mbox{\ for all\ }i\in\a\} .\] Let $\Bb$ denote the cylindric set
algebra with base set $U$ and generated by $g$. In $\Bb$, the
element $g$ is an atom, this can be seen by using permutations of
$U$ exactly as in the proof of Lemma~\ref{gen-lem}. Now, $g$ is
below all the diversity elements $-d_{ij}$, $i<j<\a$, so we can
split it into two disjoint parts $g', g''$ obtaining the algebra
$\Aa\in\CA_\a$ defined as follows.
Let $\langle A,+,-\rangle$ be the Boolean algebra which contains
$\langle B,+,-\rangle$ as a subalgebra, in which $g'$ and $g''$ are
disjoint nonzero elements such that $g=g'+g''$ and which is
generated by $B\cup\{ g'\}$. Then the elements of $A$ are
\[ A = \{ b+h : b\in B\mbox{\ and\ }h\in\{0,g',g''\}\}.\]
The  cylindric operations are defined in $\Aa$ so that $\Bb$ is a
subalgebra of $\Aa$ and
\[ c_i(b+g') \de c_i(b+g'') \de c_i(b+g)\qquad\mbox{\ for all\ }b\in
B .\] Now,  $\Aa\in\CA_\a$ can be checked directly by checking that
the cylindric equations $(C_0) - (C_7)$  of \cite[1.1.1]{HMTI} hold
in $\Aa$, or by checking that $\Aa$ is the algebra we get from $\Bb$
by splitting $g$ in it by the method in \cite[2.6.12]{HMTI}.

We show that $\Aa\not\models e_{01}(g,g')$. Indeed, $c_1(g\cdot
g'\cdot c_0(g-g'))=c_1(g'\cdot c_0(g''))=c_1(g'\cdot c_0g)=c_1g$
while $c_0(c_1g-d_{01})=c_0(V_0\times (U-V_0)\times V_2\times\dots)$
which does not contain $c_1g=V_0\times U\times V_2\times\dots$.

Next we show that $\Aa\models e_{ij}(x,y)$ when $i\ne 0$, i.e., we
show
\[ \Aa\models c_i(x\cdot y\cdot c_j(x-y))\le c_i(c_jx-d_{ij})
.\] Let $x,y\in B$ be arbitrary. Then $x\cdot y, x-y$ are of form
$a+h$, $b+k$ with $a,b\in B$,\ $h,k\in\{ 0,g',g''\}$ and $a,b,g$
pairwise disjoint as well as $h,k$ disjoint, by our construction of
$\Aa$. Since negation $-$ occurs in the equation only in form of
$-d_{ij}$, the terms at the two sides of the equation are additive,
and since $a,b\in B$ and $\Bb\subseteq\Aa$ is representable, we have
that the equation is true for $a,b$.  So, if both $h$ and $k$ are
$0$, then we are done. Assume therefore that $h+k\ne 0$. Then we get
a bigger term on the lhs of the inequality if we replace $h,k$ with
$g,g$ respectively. We get then $c_i(x\cdot y\cdot c_j(x-y)) =
c_i((a+g)\cdot c_j(b+g))=c_i(a\cdot c_j(b+g))+c_i(g\cdot c_j(b+g))$.
On the other side of the inequality we have $c_i(c_jx-d_{ij}) =
c_i(c_j(a+b+h+k)-d_{ij}) = c_i(c_j(a+b+g)-d_{ij})$. (We used
$c_j(h+k)=c_jg$ in the last step.) This is now an equation
concerning the representable algebra $\Bb$ since all the elements
occurring are in $\Bb$. Now, $c_i(a\cdot c_j(b+g))\le
c_i(c_j(a+b+g)-d_{ij})$ since this is an instance of Henkin's
equation by $a$ and $b+g$ being disjoint. We only have to show
$c_i(g\cdot c_j(b+g))\le c_i(c_j(a+b+g)-d_{ij}$. The inequality
\begin{equation}\label{g-eq}
c_i(g)\le c_i(c_jg-d_{ij})
\end{equation}
holds because $V_i$ has at least two elements: $c_i(c_jg-d_{ij}) =
c_i\{ s\in V : (\forall k\ne j)s_k\in V_k\mbox{\ and\ }s_j\ne s_i\}
= c_ig$.  Then
\begin{align*}
&c_i(g\cdot c_j(b+g))&&=&&\quad&&\mbox{by $g$ being an atom}\\
&c_i(g)&&\le&&\quad&&\mbox{by \eqref{g-eq}}\\
&c_i(c_jg-d_{ij})&&\le&&\quad&&\mbox{by monotony of the terms involved}\\
&c_i(c_j(a+b+g)-d_{ij}) .&&\quad&& &&\mbox{}
\end{align*}
To finish the proof of Theorem~\ref{ca-thm}, let $E\de\{ e_{i1} :
i\in\a, i\ne 1\}$ and for all $H\subseteq E$ let $\V_H$ be the
subvariety of $\CA_\a$ axiomatized by $H$. Then
$\V_H\supseteq\RCA_\a$ for all $H\subseteq E$, by $\RCA_\a\models
E$. Also, $\V_H\ne\V_G$ for distinct $H,G\subseteq E$ since if, say,
$e_{i1}\in G-H$ then $\Rd^{\rho}\Aa\in(\V_H-\V_G)$ whenever $\rho$
is a permutation of $\a$ with $\rho(i)=0$. \QED

\section{Counterpoint: classes with only continuum many varieties}\label{counter-sec}


Let us call a cylindric algebra $\Aa$ \emph{symmetric} iff
$\Aa\models e$ implies $\Aa\models\rho(e)$ for all permutations
$\rho$ of $\a$, where $\rho(e)$ denotes the equation we get from $e$
by systematically replacing each index $i\in\a$ in it with
$\rho(i)\in\a$.
The proofs of the previous theorems were based on the existence of
non-symmetric algebras. We will show that, surprisingly, many
$\CA_\a$s, almost all in some sense, are symmetric. In particular,
all dimension-complemented, all diagonal cylindric algebras, and
more generally, all the algebras occurring in
\cite[Thm.2.6.50(i)-(iii)]{HMTI} are symmetric. Clearly, symmetric
algebras can generate at most continuum many varieties since their
equational theories are determined by equations written in the first
$\omega$ indices. In section~\ref{fewvar-ssec} we show that this
maximal possible number $2^\omega$ is indeed achieved using only a
small subclass of symmetric $\CA_\a$s: locally finite dimensional
regular cylindric set algebras with infinite bases generate indeed
continuum many varieties, for all infinite $\a$.

Thus, $\RCA_\a$ has $2^\a$ subvarieties, but locally finite
dimensional ones generate only $2^\omega$ many. What is the property
that the $\Lf_\a$-generated varieties have but not all of the
subvarieties have? Clearly, being symmetric is such a distinguishing
property. (We call a variety symmetric iff it is generated by
symmetric algebras.) However, being symmetric does not characterize
the $\Lf_\a$-generated subvarieties: we will show that there is a
symmetric subvariety of $\RCA_\a$ that is not generated by a
subclass of $\Lf_\a$.
In section~\ref{ind-ssec} we introduce the notion of
\emph{inductive} algebras and inductive varieties and we prove that
this property characterizes the subvarieties generated by $\Lf_\a$s,
the property of being inductive singles out the $2^\omega$ many
$\Lf_\a$-generated subvarieties among all the $2^\a$ many
subvarieties of $\RCA_\a$.
By this we also get a simple characterization, and recursive
enumeration, of the equational theory of $\RCA_\a$, much simpler
than either one of the three enumerations presented in
\cite[pp.112-119]{HMTII}. These results contribute to solving
\cite[Problem 4.1]{HMTII} which is asking for a simple equational
basis for $\RCA_\a$.

Being inductive is a nice property: inductive algebras are all
representable, they are symmetric, and their equational theories
coincide with the one of an $\Lf_\a$. We show that there are more
inductive algebras than the widest class dealt with in
\cite[2.6.50]{HMTI}. This provides us with a new representation
theorem for $\CA_\a$s. All this strengthen, extend and improve
\cite[2.6.50]{HMTI}, whose significance was discussed in the
introduction to the present paper. The notion of being inductive can
be described by a set of $\Delta_2$ first order logic formulas.
Since inductive algebras are all symmetric and we have constructed a
nonsymmetric representable algebra in section~\ref{const-sec} here,
we get a $\Delta_2$-formula distinguishing $\Lf_\a$ and $\RCA_\a$.

\subsection{Endo-dimension-complemented algebras are
symmetric}\label{sym-ssec}
Let $\Lf_\a$, $\Dc_\a$, and $\Di_\a$ denote the classes of all
locally finite dimensional, dimension-complemented, and diagonal
$\CA$'s, respectively. Let us call the elements of the wider class
introduced in (iii) of \cite[2.6.50]{HMTI}
\emph{endo-dimension-complemented} (\emph{endo-dc} in short): an
algebra $\Aa\in\CA_\a$ is called endo-dc if for each finite
$\Gamma\subseteq\a$ and each nonzero $x\in A$ there are a
$\kappa\in\a-\Gamma$ and an endomorphism $h$ of the
\emph{$\Gamma$-reduct} $\Rd_{\Gamma}\Aa\de\langle A,+,-,c_i,
d_{ij}\rangle_{i,j\in\Gamma}$ of $\Aa$ such that $h(x)\ne 0$ and
each element of the range of $h$ is $\kappa$-closed, i.e., $c_\kappa
h(a)=h(a)$ for all $a\in A$. Let $\Edc_\a$  denote the class of all
endo-dimension-complemented $\CA_\a$s. It is proved in
\cite[2.6.50]{HMTI} that $\Lf_\a\subset\Dc_\a\subset\Di_\a\subset
\Edc_\a\subseteq\RCA_\a$ and it is asked as \cite[Problem
2.13]{HMTI} whether the last inclusion is proper or not. We are
going to show that this inclusion is proper: the algebra constructed
in section~\ref{const-sec} here is representable but not endo-dc.
More specifically, we will show that each endo-dc algebra is
symmetric, which implies $\Edc_\a\ne\RCA_\a$ since our witness
algebra $\Aa$ in the proof of Theorem~\ref{manyvar-thm} was designed
to be non-symmetric but it is representable. We also show that
$\RCA_\a$ is close to $\Edc_\a$ in the sense that $\RCA_\a$ is the
closure of $\Edc_\a$ under taking subalgebras. On the other hand, to
indicate the distance between $\Edc_\a$ and $\RCA_\a$ we show that
the class $\Sy_\a\cap\RCA_\a$ of symmetric representable algebras
lies strictly in between $\Edc_\a$ and $\RCA_a$, i.e.,
$\Edc_\a\subset\Sy_\a\cap\RCA_\a\subset\RCA_\a$.

\begin{thm}\label{sym-thm} Each endo-dc algebra is symmetric.
\end{thm}

\noindent{\bf Proof.}
Because the notion of a symmetric algebra involves renaming indices
of operations, in this and the coming proofs we will often deal with
renaming operations in equations and algebras. Therefore we begin
the proof with introducing notation for these. We will use these
notation, except for $\ind(e)$, only in proofs.

If $\tau$ is a term in the language of $\CA_\a$ and $\rho:\a\to\a$
then $\rho(\tau)$, the term we get from $\tau$ by renaming the
indices occurring in it according to $\rho$, is defined by induction
as $\rho(d_{ij})\de d_{\rho(i)\rho(j)}$, $\rho(c_i\sigma)\de
c_{\rho(i)}\rho(\sigma)$ and $\rho(x)\de x$ if $x$ is a variable,
$\rho(\sigma+\delta)\de\rho(\sigma)+\rho(\delta)$,
$\rho(-\sigma)\de-\rho(\sigma)$. If $e$ is an equation of form
$\tau=\sigma$ then $\rho(e)$ is $\rho(\tau)=\rho(\sigma)$.

$\ind(\tau)$ denotes the set of indices occurring in $\tau$, this is
defined by induction as follows. $\ind(d_{ij})\de\{ i,j\}$,
$\ind(c_i\tau)\de\{ i\}\cup\ind(\tau)$, and $\ind(x)\de\emptyset$,
$\ind(\tau+\sigma)\de\ind(\tau)\cup\ind(\sigma)$,
$\ind(-\tau)\de\ind(\tau)$. If $e$ is an equation of form
$\tau=\sigma$ then $\ind(e)\de\ind(\tau)\cup\ind(\sigma)$.

Assume $\Aa\in\CA_\a$, $\Gamma$ is any set and $\rho:\Gamma\to\a$ is
a one-to-one function. Then $\Rd^{\rho}\Aa$ denotes an algebra whose
signature is that of $\CA_{\Gamma}$, whose Boolean reduct $\langle
A,+,-\rangle$ is the same as that of $\Aa$, whose operation denoted
by $c_i$ for $i\in\Gamma$ is the operation of $\Aa$ denoted by
$c_{\rho(i)}$, and similarly for the diagonals, $d_{ij}$ of
$\Rd^{\rho}\Aa$ is the same as $d_{\rho(i)\rho(j)}$ of $\Aa$. In
symbols, \[ \Rd^{\rho}\Aa\de\langle
A,+,-,c_{\rho(i)},d_{\rho(i)\rho(j)} : i,j\in\Gamma\rangle .\] It is
not difficult to check that $\Rd^{\rho}\Aa\in\CA_{\Gamma}$ and
$\Rd^{\rho}\Aa\models e$ iff $\Aa\models \rho(e)$, for any equation
$e$. This algebra is called a generalized reduct of $\Aa$ and it is
introduced in \cite[2.6.1]{HMTI}.\smallskip

We begin the proof of Theorem~\ref{sym-thm}.
Assume $\Aa\in\CA_\a$ is endo-dc, we want to show that it is
symmetric. This means showing that $\Aa\models e$ implies
$\Aa\models\rho(e)$ for all equations $e$ and permutations $\rho$ of
$\a$. For this, it is enough to prove
\begin{equation}\Aa\not\models e\quad\text{implies}\quad
\Aa\not\models\rho(e), \qquad\text{ for all }e\text{ and }\rho,
\end{equation}
since each equation $e$ is of form $\rho(e')$ and
$\rho^{-1}\rho(e')=e'$. Assume $\Aa\not\models e$. We may assume
that $e$ is of form $\tau=0$ for some $\tau$.
Let $\Gamma\de\ind(\tau)$ and let $\rho:\Gamma\to\Delta$ be a
bijection. We have $\tau(a)\ne 0$ for some $a\in A$ by
$\Aa\not\models e$, and we want to show that $\rho(\tau)(b)\ne 0$
for some $b\in A$. (In fact, $\tau$ may have more than one variable,
so we should use a sequence $\overline{a}$ in place of $a\in A$. For
simplicity, we write out the present proof for the case when $\tau$
contains one variable.)

We aim for getting a homomorphism $\Rd_\Gamma\Aa\to\Rd^{\rho}\Aa$
which takes $\tau(a)$ to a nonzero element.
The idea of the proof is as follows. Assume $\Delta=\{
k_1,\dots,k_n\}$ is disjoint from $\Gamma=\{ i_1,\dots,i_n\}$. Then
the substitution operation $x\mapsto \s^{i_1}_{k_1}\dots
\s^{i_n}_{k_n}(x)$ is such a homomorphism, but only on the
$\Delta$-closed elements $x$, i.e., when $x=c_{(\Delta)}x\de
c_{k_1}\dots c_{k_n}x$. There are two obstacles to deal with:
$\Delta$ may not be disjoint from $\Gamma$, and $\tau(a)$ may not be
$\Delta$-closed. We deal with the first obstacle by finding $J$
which is disjoint both from $\Gamma$ and $\Delta$, and finding
desired homomorphisms from $\Gamma$ to $J$ and then from $J$ to
$\Delta$. We deal with the second obstacle by using the condition
$\Aa\in\Edc_\a$ for finding a homomorphism from $\Gamma$ to $\Gamma$
which takes $\tau(a)$ to a $J$-closed non-zero element. We begin now
to elaborate the just outlined idea.

By \cite[(2), p.416]{HMTI}, $\Aa\in\Edc_\a$ implies that there is
$J\subseteq\a-(\Gamma\cup \Delta)$ with $|J|=|\Gamma|$ and there is
a homomorphism $h:\Rd_{\Gamma}\Aa\to \Rd_{\Gamma}\Aa$  such that
$h(\tau(a))\ne 0$ and  $h(x)=c_{(J)}h(x)$ for all $x\in A$. By $h$
being a homomorphism on $\Rd_{\Gamma}\Aa$ and
$\ind(\tau)\subseteq\Gamma$ we have $\tau(h(a))=h(\tau(a))\ne 0$.

Now that $h$ provided us with $J$-closed elements, we can use the
usual substitution operations $\s^i_j$ to get the homomorphism we
seek for, as follows. Let $c_j^*\Aa$ denote the algebra whose
elements are the $c_j$-closed elements of $\Aa$ and whose operations
are those of $\Aa$ except $c_j, d_{jk}, d_{kj}$ for $k\in\a$.  This
is indeed an algebra, it is $\Nr_{(\alpha-\{j\})}\Aa$ in the
terminology of \cite{HMTI}, but we will use the shorter notation
$c_j^*\Aa$ in the present proof. We will use $c_{(J)}^*\Aa$ for the
analogous algebra (where $J\subset\a$). Let $[i/j]$ denote the
function that takes $i$ to $j$ and takes $k$ to $k$ for all
$k\in(\a-\{ i,j\})$. Then $\Rd^{[i/j]}c_i^*\Aa$ is the algebra
$c_i^*\Aa$ except that we rename the operations $c_j, d_{jk},
d_{kj}$ (of $c_i^*\Aa$) as $c_i, d_{ik}, d_{ki}$, respectively. Thus
the similarity types of $c_j^*\Aa$ and $\Rd^{[i/j]}c_i^*\Aa$ are
equal.
We are going to show, by using \cite[sec.\ 1.5]{HMTI}, that
\begin{equation}\label{subs-eq}
\s^i_j: c_j^*\Aa \to \Rd^{[i/j]}c_i^*\Aa\quad\text{ is an
isomorphism.}
\end{equation}
Indeed, $\s^i_j$ is a Boolean homomorphism by 1.5.3, it is a
homomorphism for $c_k, d_{km}$ for $k,m\in\a-\{i,j\}$ by 1.5.8(ii),
1.5.4(ii), and it takes $d_{ik}, d_{ki}$ to $d_{jk}, d_{kj}$ by
1.5.4(i). For the next two steps we need to use that we are mapping
$c_j$-closed elements. $\s^j_i$ is the inverse of $\s^i_j$ on
$c_j$-closed elements because $\s^j_i\s^i_jc_jx=c_jx$ by 1.5.10(i),
1.5.8(i). $\s^i_j$ takes the operation $c_i$ on $c_j$-closed
elements to $c_j$ because 
$c_j\s^i_ja = c_i\s^j_ia = c_i\s^j_ic_ja = c_ic_ja = c_ia =
\s^i_jc_ia$, by 1.5.8(i), 1.5.9(i). We are done with proving
\eqref{subs-eq}.

Recall that $J\subseteq\a-(\Gamma\cup\Delta)$ and $|J|=|\Gamma|$.
Let $i_1,\dots,i_n$ and $j_1,\dots,j_n$ be repetition-free
enumerations of $\Gamma$ and $J$, respectively. Let $\eta:\Gamma\to
J$ be such that $\eta(i_1)=j_1, \dots, \eta(i_n)=j_n$. Define
\[ \s(\eta) \de \s^{i_1}_{j_1}\dots \s^{i_n}_{j_n}  .\]
By using \eqref{subs-eq} successively, we get
\begin{equation}\label{subse-eq}
\s(\eta): \Rd_\Gamma c_{(J)}^*\Aa \to
\Rd^{\eta}c_{(\Gamma)}^*\Aa\quad\text{ is an isomorphism.}
\end{equation}
By letting $k_\ell\de\rho(i_\ell)$ and $\xi(j_\ell)\de k_\ell$ for
$1\le\ell\le n$ we get that $k_1,\dots, k_n$ is a repetition-free
enumeration of $\Delta$, and $\rho = \xi\circ\eta$. By repeating the
process leading to \eqref{subse-eq} we get
\begin{equation}\label{subsx-eq}
\s(\xi): \Rd^{\eta}c_{(\Gamma)}^*\Aa \to
\Rd^{\rho}c_{(\Delta)}^*\Aa\quad\text{ is an isomorphism.}
\end{equation}
Putting these two isomorphisms together we get
\begin{equation}\label{subsr-eq}
\s(\rho) : \Rd_\Gamma c_{(J)}^*\Aa \to
\Rd^{\rho}c_{(\Delta)}^*\Aa\quad\text{ is an isomorphism.}
\end{equation}
Let $g\de \s(\rho)\circ h$, then $g(\tau(a))=\s(\rho)h(a)\ne 0$ by
$h(a)\ne 0$, so
\begin{equation}\label{end-eq}
g:\Rd_\Gamma \Aa\to \Rd^{\rho}c_{(\Delta)}^*\Aa\quad\text{ is a
homomorphism with}\ \ g(\tau(a))\ne 0.
\end{equation}
Now,
$\rho(\tau)$ in $\Aa$ is the same as $\tau$ in $\Rd^\rho\Aa$, by
definition. Therefore, $\rho(\tau)(ga)$ in $\Aa$ is the same as
$\tau(ga)$ in $\Rd^\rho\Aa$, which is the same as $g(\tau(a))$ which
is nonzero by \eqref{end-eq} and $\tau(a)\ne 0$. We are done with
showing that $\Aa$ is symmetric. \QED

\begin{lem}\label{full-lem}
Each full cylindric set algebra with unit a disjoint union of weak
spaces is endo-dc.
\end{lem}

\noindent{\bf Proof.} The proof in \cite[2.6.51, p.417]{HMTI} for
showing ``(iii) does not imply (ii)" in fact proves the present
Lemma~\ref{full-lem}.\QED

\begin{thm}\label{edc-thm}\text{\rm (Solution of \cite[Problem
2.13]{HMTI})} There is an $\RCA_\a$ which is not endo-dc, but each
$\RCA_\a$ can be embedded into an endo-dc algebra. In symbols:
$\Edc_\a\subset\Sub\Edc_\a=\RCA_\a$.
\end{thm}

\noindent{\bf Proof.} The algebra we based the proof of
Theorem~\ref{manyvar-thm} on is not symmetric, hence not endo-dc by
Theorem~\ref{sym-thm}. Clearly, it is representable. This shows
$\Edc_\a\ne\RCA_\a$. $\Edc_\a\subseteq\RCA_\a$ is proved as
(iii)$\Rightarrow$(iv) in \cite[Thm.2.6.50]{HMTI}.
$\RCA_\a=\Sub\Edc_\a$ follows from Lemma~\ref{full-lem} immediately,
since each representable algebra is embeddable into a full one.\QED

\begin{rem}\label{instance-rem}
\rm{In the proof above, we used Theorem~\ref{sym-thm} to show that
the algebra $\Aa$ we used in the proof of Theorem~\ref{manyvar-thm}
is not endo-dc. A concrete $\Gamma\subseteq\a$ and nonzero $a\in A$
for which there are no $\kappa\in\a$ and endomorphism $h$ with the
required properties are $\{ 0,1,2\}$ and $g$. Indeed, let $\tau\de
x-c_{(3)}(\beta+\gamma+\o)$, see \eqref{cutt-eq}. Then $e(x)$ fails
iff $\tau(x)\ne 0$, by \eqref{cutt-eq}. Hence, $\tau(g)\ne 0$ but
$\tau(c_\kappa x)=0$ for all $\kappa$ by
Remark~\ref{constr-rem}(iii), and this implies that there is no
endomorphism $h$ of $\Rd_{\{0,1,2\}}\Aa$ with range inside
$c_\kappa^*A$ and $h(g)\ne 0$.}
\end{rem}

\begin{thm}\label{sym-lem}
Not all symmetric algebras are representable, and not all
representable algebras are symmetric. In symbols,\\
\[\Sy_\a\cap\RCA_\a\subset\Sy_\a\quad\text{ and }\quad
\Sy_\a\cap\RCA_\a\subset\RCA_\a. \]
\end{thm}

\noindent{\bf Proof.} To exhibit a symmetric algebra that is
nonrepresentable, take any nonrepresentable $\Aa\in\CA_\a$, we
``turn" it symmetric. Indeed, let
\[ \Bb \de \prod\langle \Rd^\rho\Aa : \rho\mbox{\ is a permutation
of\ }\a \rangle .\] That $\Bb$ is symmetric can be seen by
\begin{align*}
&\Bb\models e&&\mbox{iff} &&\quad&&\mbox{by the definition of $\Bb$}\\
&\Rd^\rho\Aa\models e\mbox{\ for all $\rho$}&&\mbox{iff}&&\quad&&\mbox{by the definition of $\Rd^\rho$}\\
&\Aa\models\rho(e)\mbox{\ for all $\rho$}&&\mbox{iff}&&\quad&&\mbox{by the nature of permutations}\\
&\Aa\models\rho(\eta(e))\mbox{\ for all $\rho,\eta$}&&\mbox{iff}&&\quad&&\mbox{by previous step}\\
&\Rd^\rho\Aa\models\eta(e)\mbox{\ for all $\rho$}&&\mbox{iff}&&\quad&&\mbox{by first step}\\
&\Bb\models\eta(e) .&&\quad&& &&\mbox{}
\end{align*}
That $\Bb$ is not representable follows from the facts that
$\Aa\notin\RCA_\a$ is a homomorphic image of $\Bb$ (as
$\Aa=\Rd^\rho\Aa$ with $\rho$ being the identity permutation of
$\a$) and $\RCA_\a$, being a variety (\cite[3.1.103]{HMTII}), is
closed under homomorphic images.

The algebra used in the proof of Theorem~\ref{manyvar-thm} is
representable and non-symmetric, this proves the second part of the
theorem, i.e., $\RCA_\a\cap\Sy_\a\subset\RCA_\a$.\QED

\subsection{Polyadic algebras are symmetric}\label{poly-ssec}

We have seen in the proof of Theorem~\ref{sym-thm} that substitution
operations are useful in proving an algebra be symmetric. In fact,
the proof of Theorem~\ref{manyvar-thm} hinges over the fact that the
polyadic substitution operations $\p_{ij}$ are not expressible in
the witness algebra $\Aa$. In this section we very briefly talk
about Halmos' polyadic algebras. We show that $\a$-dimensional
quasi-polyadic equality algebras indeed have only $2^\omega$ many
subvarieties, since  all their members are symmetric (in an
appropriate sense). We then state some of the corollaries of our
construction that concern polyadic algebras.

Polyadic equality algebras ($\PEA_\a$s) were introduced by Paul
Halmos \cite{Hal62}, they are basically cylindric algebras endowed
with unary substitution operations $\s_\rho$ for $\rho:\a\to\a$. In
the set algebras with unit ${}^\a U$ these are interpreted as
\[ S_\rho(X)\de\{ s\in{}^\a U : \rho\circ s \in X\}.\]
Quasi-polyadic equality algebras were also defined by Halmos in
\cite{Hal62}, they retain only those substitutions where $\rho$ is
finite. Let $\QPEA_\a$ denote their class, for precise definition
see, e.g., \cite[p.266, item 9]{HMTII} or \cite{ST}.

\begin{thm}\label{poly-thm}
$\QPEA_\a$ has exactly $2^\omega$ many subvarieties.
\end{thm}

\noindent{\bf Proof.} The idea of the proof is to show that each
$\QPEA_\a$ is symmetric in the sense analogous to the notion used in
$\CA_\a$. However, the indices of the $\QPEA_\a$-operations have
some structure, it is not so clear how we are to change the indices
in an equation systematically/uniformly. (For more on this see the
introduction of \cite{ST}.) Therefore, we will use the more
index-friendly version $\FPEA_\a$ of $\QPEA_\a$ defined in
\cite{ST}. Since the two varieties are term-definitionally
equivalent, proved as \cite[Thm.1(ii)]{ST}, it is enough to show
that $\FPEA_\a$ has only $2^\omega$ subvarieties.

We are going to show that each element of $\FPEA_\a$ is symmetric in
the very analogous sense to $\CA_\a$, this will prove that
$\QPEA_\a$ has at most continuum many subvarieties (since each
equation is equivalent to one which uses indices from $\omega$
only). That $\QPEA_\a$ has indeed continuum many varieties can be
seen by repeating the proof of \cite[Thm.4.1.24]{HMTII} for
$\QPEA_\a$.

The extra-cylindric operations in an $\FPEA_\a$ are denoted as
$\p_{ij}$ for $i,j\in\a$. The operation $\p_{ij}$ stands for
$\s_{\rho}$ where $\rho$ is $[i,j]$, the latter being the
permutation of $\a$ that interchanges $i$ and $j$ and leaves all the
other elements fixed.  Now, the definitions of $\rho(\tau)$ and
$\ind(\tau)$ for $\FPEA_\a$-terms $\tau$ can easily be extended from
$\CA_\a$.
We will show that each $\Aa\in\FPEA_\a$ is symmetric in the sense
that
\begin{equation*}
\Aa\models e\mbox{\quad iff\quad}\Aa\models\rho(e),\qquad\mbox{ for
all permutations $\rho$ of $\a$.}
\end{equation*}
Indeed, let $\Aa\in\FPEA_\a$ and let $\tau$ be a term in the
language of $\FPEA_\a$, let $\rho$ be a permutation of $\a$. Then
$\ind(\tau)$ is finite, so we may assume that $\rho$ is finite, too.
Each finite permutation is a composition of transpositions $[i,j]$,
so we may assume that $\rho$ is indeed a transposition $[i,j]$ with
$i\ne j$. In the sequel we will write $\tau(\bar{x})$ and
$\tau(\p_{ij}\bar{x})$ for $\tau(x_1,\dots,x_n)$ and
$\tau(\p_{ij}x_1,\dots,\p_{ij}x_n)$. The following can be proved by
induction on $\tau$:
\begin{equation}\label{rho-eq} \FPEA_\a\models\ \
\p_{ij}\tau(\bar{x})=\rho(\tau(\p_{ij}\bar{x}))\end{equation} with
the use of the following equations that can be proved to hold in
$\FPEA_\a$:
\begin{description}
\item[]
$\p_{ij}(x+y)=\p_{ij}x+\p_{ij}y,\quad \p_{ij}(-x)=-\p_{ij}x, \quad$
\item[]
$\p_{ij}\p_{ij}x=x,\quad \p_{ij}x=\p_{ji}x,\quad$
\item[]
$\p_{ij}c_kx=c_{k'}\p_{ij}x,\quad\p_{ij}d_{kl}=d_{k'l'},\quad\p_{ij}\p_{kl}x=\p_{k'l'}\p_{ij}x,\quad$
\end{description}
where $k'=\rho(k)$ and $l'=\rho(l)$.  Now, let  $e$ be any equation,
we may assume that it is of form $\tau(\bar{x})=1$.
\begin{align*}
&\Aa\models e&&\mbox{iff} &&\quad&&\mbox{by $e$ being $\tau=1$}\\
&\Aa\models \tau=1&&\mbox{implies} &&\quad&&\mbox{by $\p_{ij}1=1$}\\
&\Aa\models \p_{ij}\tau=1&&\mbox{iff}&&\quad&&\mbox{by \eqref{rho-eq}}\\
&\Aa\models\rho(\tau(\p_{ij}\bar{x}))=1&&\mbox{implies}&&\quad&&\mbox{by $\p_{ij}\p_{ij}x=x$}\\
&\Aa\models\rho(\tau)=1&&\mbox{iff}&&\quad&&\mbox{by $e$ being $\tau=1$}\\
&\Aa\models\rho(e) .&&\quad&& &&\mbox{}
\end{align*}\QED

It is proved in \cite[5.4.18]{HMTII} that the cylindric reducts of
$\PEA_\a$s are all representable, in symbols
$\Rrd_{ca}\PEA_\a\subseteq\RCA_\a$. Our results imply that this
inclusion is a strict one. Further, the cylindric reducts of
(quasi-)polyadic (equality)-algebras are not closed under
subalgebras.

\begin{cor}\label{cared-cor}
Not every representable cylindric algebra is the cylindric reduct of
a polyadic equality algebra, hence the class of the latter is not
closed under subalgebras. Formally:
\[\Rrd_{ca}\PEA_\a\subset\RCA_\a=\Sub\Rrd_{ca}\PEA_\a.\]
Further,  $\Rrd_{ca}\QPEA_\a\subset\Sub\Rrd_{ca}\QPEA_\a$.
\end{cor}

\noindent{\bf Proof.} It follows from the proof of
Theorem~\ref{poly-thm} that the cylindric reduct of any
quasi-polyadic equality algebra is symmetric. We have seen in
Theorem~\ref{sym-lem} that not all representable algebras are
symmetric. Take a non-symmetric $\RCA_\a$, it is not in
$\Rrd_{ca}\QPEA_\a$, hence it is not in $\Rrd_{ca}\PEA_\a$, either.
Since all full cylindric set algebras are reducts of $\PEA_\a$, our
non-symmetric $\RCA_\a$ is in $\Sub\Rrd_{ca}\PEA_\a$. \QED

\subsection{Inductive algebras}\label{ind-ssec}

Let us call a cylindric algebra $\Aa$ \emph{inductive} iff
$\Aa\models e(c_ix_1,\dots,c_ix_n)$ implies $\Aa\models
e(x_1,\dots,x_n)$ whenever $e$ is an equation and $i$ does not occur
as an index of an operation in $e$. Let $\Ret_\a$ denote the class
of all inductive $\CA_\a$s. While $\Aa\models e$ implies $\Aa\models
e(c_ix)$ always holds, the converse of this would be thought to hold
only in rather special cases, if at all. We are going to show that,
on the contrary, there is a great variety of inductive algebras:
each endo-dc algebra is inductive and we have already seen that
there is a great variety of endo-dc algebras. There are even more
inductive algebras than endo-dc algebras: $\Edc_\a\subset\Ret_\a$.
We then prove that each inductive algebra is representable and
symmetric (but the converse does not hold). Thus, we refine the
chain
$\Lf_\a\subset\Dc_\a\subset\Di_\a\subset\Edc_\a\subset\Sy\cap\RCA_\a\subset\RCA_\a$
with inserting a new class into it: $\Edc_\a\subset\Ret_\a\subset
\Sy_\a\cap\RCA_\a$. This is also a new representation theorem, a
sharpening of \cite[2.6.50]{HMTI}, since in the chain presented in
\cite[2.6.50]{HMTI} the widest representable class was $\Edc_\a$.
The new class $\Ret_\a$ has an additional significance, namely an
algebra is inductive iff it is equationally indistinguishable from
an $\Lf_\a$. So, inductive algebras are in intimate connection with
$\Lf_\a$. This will give us a specific $\Delta_2$ formula
distinguishing $\Lf_\a$ and $\RCA_\a$.

\goodbreak

\begin{thm}\label{ret-thm}\hspace{10pt}
\begin{description}
\item{(i)}
Each endo-dc algebra is inductive, and each inductive algebra is
symmetric and representable but the converses of these statements do not hold, i.e., 
\[\Edc_\a\subset\Ret_\a\subset\Sy_\a\cap\RCA_\a.\]
\item{(ii)}
An algebra is inductive iff there is an $\Lf_\a$ with the same
equational theory, i.e.,
\[ \Aa\text{ is inductive }\quad\text{ iff }\quad \Eq(\Aa)=\Eq(\Bb)\text{ for
some }\Bb\in\Lf_\a .\]
\end{description}
\end{thm}

\noindent{\bf Proof.} First we prove part of (i), namely we prove
$\Edc_\a\subseteq\Ret_\a$. This follows almost directly from the
definitions and from Theorem~\ref{sym-thm}. Let $\Aa\in\Edc_\a$ and
let $e(x_1,...,x_n)$ be an equation, $i\in\a$ such that $i$ does not
occur in $e$. In the sequel we will write $e(\bar{x})$ and
$e(c_i\bar{x})$ in place of $e(x_1,\dots,x_n)$ and
$e(c_ix_1,\dots,c_ix_n)$, respectively. We want to show $\Aa\models
e(c_i\bar{x})$ implies $\Aa\models e(\bar{x})$. To this end, we
assume $\Aa\not\models e(\bar{x})$ and we show that $\Aa\not\models
e(c_i\bar{x})$. Let $a_1,\dots,a_n\in A$ be such that
$\Aa\not\models e(\bar{a})$. We may assume that $e$ is of form
$\tau=0$, so we have $\tau(\bar{a})\ne 0$ in $\Aa$. Let
$\Gamma\de\ind(\tau)$. By $\Aa$ being endo-dc, there are a
homomorphism $h:\Rd_\Gamma\Aa\to\Rd_\Gamma\Aa$ and a
$\kappa\in\a-\Gamma$ such that $h(\tau(\bar{a}))\ne 0$ and
$h(b)=c_\kappa h(b)$ for all $b\in A$. Now,
$h(\tau(\bar{a}))=\tau(h(\bar{a}))=\tau(h(a_1),\dots,h(a_n))$ by $h$
being a homomorphism wrt.\ the operations occurring in $\tau$. By
$h(a_1)=c_\kappa h(a_1),\dots , h(a_n)=c_\kappa h(a_n)$ we then have
$\tau(c_\kappa h(\bar{a}))\ne 0$ in $\Aa$. This means that
$\Aa\not\models e(c_\kappa \bar{x})$. Since $\Aa$ is symmetric by
Theorem~\ref{sym-thm} and $\kappa, i\notin\ind(e)$, we get that
$\Aa\not\models e(c_i\bar{x})$ as was desired.

Next we prove (ii). For proving the ``only-if\," part, let $\Aa$ be
inductive, we will show that it is equationally indistinguishable
from an $\Lf_\a$. Let $\Cc$ be an elementary $\a$-saturated
extension of $\Aa$, and let $\Bb$ be the greatest locally finite
dimensional subalgebra of $\Cc$. (This exists by
\cite[2.1.5(ii)]{HMTI}.) We are going to show that $\Aa$ and $\Bb$
are equationally indistinguishable. Let $e$ be an equation. If
$\Aa\models e$ then $\Cc\models e$ because $\Cc$ is an elementary
extension of $\Aa$, and thus $\Bb\models e$ because $\Bb$ is a
subalgebra of $\Cc$. Assume now $\Aa\not\models e$. Let $\Delta\de\{
i_1,...,i_n\}$ be disjoint from $\ind(e)$ with $i_1,...,i_n$ being
all distinct. Then $\Aa\not\models e(c_{i_1}\bar{x})$ since $\Aa$ is
inductive and $\Aa\not\models e$. But then $\Aa\not\models
e(c_{i_1}c_{i_2}\bar{x})$ because
$i_2\notin\ind(e(c_{i_1}\bar{x}))$, and so on, showing that
$\Aa\not\models e(c_{(\Delta)}\bar{x})$. Let
\[\Sigma(\bar{x})\de\{ \neg e(\bar{x}), c_i\bar{x}=\bar{x} : i\in\a-\ind(e) \} .\] Each
finite subset of $\Sigma$ is satisfiable in $\Cc$ by $\Aa\not\models
e(c_{(\Delta)}\bar{x})$ for all finite $\Delta\subseteq\a-\ind(e)$.
By $\Cc$ being $\a$-saturated, this implies that there are
$b_1,\dots, b_n\in C$ for which $\Sigma(\bar{b})$ holds in $\Cc$.
These $b_j$s are finite dimensional (by $c_i(b_j)=b_j$ for all
$i\in\a-\ind(e)$), and $\Cc\not\models e(\bar{b})$ (by $\neg
e(\bar{x})\in\Sigma(\bar{x})$). Hence $b_1,\dots, b_n\in B$ and
$\Cc\not\models e(\bar{b})$, hence $\Bb\not\models e(\bar{b})$,
i.e., $\Bb\not\models e$. This finishes the ``{only-if}\," part of
the proof of (ii). For the ``if\," part, we have to show that each
$\Bb\in\Lf_\a$ is inductive. Indeed, $\Lf_\a\subseteq\Edc_\a$ by
\cite[2.6.50]{HMTI}, and $\Edc_\a\subseteq \Ret_\a$ by that part of
(i) that we have already proved.

It remains to prove the rest of (i). We have already shown
$\Edc_\a\subseteq\Ret_\a$. To show
$\Ret_\a\subseteq\Sy_\a\cap\RCA_\a$ we use (ii), \cite[2.6.50]{HMTI}
and Theorem~\ref{sym-thm}, as follows. Let $\Aa\in\Ret_\a$. Then
$\Eq(\Aa)=\Eq(\Bb)$ for some $\Bb\in\Lf_\a$, by (ii). Now,
$\Lf_\a\subseteq\Edc_\a\subseteq\Sy_\a$ by \cite[2.6.50]{HMTI} and
Theorem~\ref{sym-thm}, $\Lf_\a\subseteq\RCA_\a$ by \cite{HMTII}.
Thus,  $\Aa\in\Sy_\a\cap\RCA_\a$. We turn to proving that the stated
inclusions are proper.

First we want to exhibit an inductive algebra that is not endo-dc.
The difference between the two notions, and this will be reflected
in the algebra $\Dd$ we exhibit, is that the notion of being
inductive talks about the equational theory of the algebra, while
the notion endo-dc talks about the inner structure of the algebra.
The algebra $\Dd$ is a direct product of the $\omega$-generated free
$\RCA_\a$-algebra $\Ff$ and another representable algebra $\Aa$. By
this, it is already ensured that $\Dd$ is inductive, as follows.
\begin{align*}
&\Dd\models e(c_i\bar{x})&&\mbox{implies} &&\quad&&\mbox{by $\Ff$ being a homomorphic image of $\Dd$}\\
&\Ff\models e(c_i\bar{x})&&\mbox{implies}&&\quad&&\mbox{by $\Ff$ being a free algebra of $\RCA_\a$}\\
&\RCA_\a\models e(c_i\bar{x})&&\mbox{implies}&&\quad&&\mbox{by $\Lf_\a\subseteq\RCA_\a$}\\
&\Lf_\a\models e(c_i\bar{x})&&\mbox{implies}&&\quad&&\mbox{by $\Lf_\a\subseteq\Ret_\a$}\\
&\Lf_\a\models e(\bar{x})&&\mbox{implies}&&\quad&&\mbox{by $\RCA_\a=\Var(\Lf_\a)$}\\
&\RCA_\a\models e(\bar{x})&&\mbox{implies}&&\quad&&\mbox{by $\Dd\in\RCA_\a$}\\
&\Dd\models e(\bar{x}) .&&\quad&& &&\mbox{}
\end{align*}
The role of the algebra $\Aa$ in the direct product is to destroy
the property endo-dc. The idea is that we split an $\a$-dimensional
atom in $\Aa$ into more parts than there are $i$-closed elements
(for some $i\in\a$) in $\Dd$, so each required endomorphism will
have to collapse all of the split parts to $0$. We begin to
elaborate this idea. Let $W$ be a set of cardinality bigger than
$|\a|$ and let $\langle W,+,z\rangle$ be any commutative group on
$W$ where $z$ is the zero-element of $+$ (i.e., $w+z=w$ for all
$w\in W$). Let $U_i\de W\times\{ i\}$, let $p\de\langle (z,i) :
i\in\a\rangle$, let $U\de \bigcup\{ U_i : i\in\a\}$ and $T\de\{
s\in{}^\a U : s_i\in U_i\ \text{for all }i\in\a\ \text{ and } |\{
i\in\a : s_i\ne p_i\}|<\omega\}$.  Let $\Bb$ be the weak cylindric
set algebra of dimension $\a$ with unit element ${}^\a U^{(p)}$ and
generated by $T$. Then $|B|=|\a|$ and $T$ is an atom in $\Bb$. We
now split $T$ in $\Bb$ into $|W|$ many parts. For all $g\in W$ let
\[ T_g \de \{ s\in T : \sum\{ w : s_i=(w,i)\ \text{ for some
}i\in\a\}=g \} .\] Then the $T_g$s ($g\in W$) form a disjoint union
of $T$ such that
\begin{equation}\label{split-eq}
c_iT_g = c_iT\quad\text{ for all }i\in\a\ \text{ and }g\in W .
\end{equation}
Let $\Aa$ be the weak set algebra with unit element ${}^\a U^{(p)}$
and generated by $T$ together with $T_g$, $g\in W$. It is not hard
to check, by using \eqref{split-eq}, that each element of $A$ is of
form \[ b+\sum\{ T_g : g\in X\}\] where $b\in B$ and $X$ is a finite
or co-finite subset of $W$. Thus, all elements of $A-B$ are
$\a$-dimensional. We now show that $\Dd=\Ff\times\Aa$ is not
endo-dc. Let's fix a $g\in W$, let $a\de\langle 0,T_g\rangle\in D$,
let $\Gamma\de \{0\}$, we want to show that there are no
endomorphism $h$ of $\Rd_{\Gamma}\Dd$ and $\kappa\in\a$ such that
$h$ takes $a$ to a nonzero element and each element of the range of
$h$ is $\kappa$-closed. For, assume the contrary, that $h$ and
$\kappa$ are as described above, we will derive a contradiction.
Since there are only $\a$ many $\kappa$-closed elements of $\Aa$,
hence of $\Dd$ by $|F|=|\a|$, and there are more than $\a$ many
split parts of $T$, the endomorphism $h$ has to take two of the
elements $\langle 0,T_w\rangle\in D$ to the same element. But these
are all disjoint from each other, so $h(\langle 0,T_w\rangle)=0$ for
some $w\in W$. But then $h(c_0\langle 0,T_w\rangle)=c_0h(\langle
0,T_w)=0$ since $h$ is a homomorphism wrt.\ $c_0$. However,
$c_0\langle 0,T_w\rangle=\langle 0,c_0T_w\rangle=\langle
0,c_0T\rangle\ge\langle 0,T\rangle\ge\langle 0,T_g\rangle$, showing
that $h(\langle 0,T_g\rangle)=h(a)=0$, and this contradicts our
assumption $h(a)\ne 0$. Thus the algebra $\Dd$ is inductive but not
endo-dc.

Finally, we exhibit a symmetric representable algebra which is not
inductive. Here, both notions refer to the equational theory of the
algebra, but they make different restrictions on it. Symmetry
requires that if an equation holds then its versions where we rename
the indices hold also, and inductivity requires that the same
equations hold for the some-cylindrification-closed elements than
for the whole algebra. Our algebra $\Aa$ that we used in the proof
of Theorem~\ref{manyvar-thm} is not symmetric, hence it is not
inductive, either, by the already proved part of (i) of the present
theorem. We will modify the algebra $\Aa$ so that it becomes
symmetric, but the above mentioned difference between the
some-cylindrification-closed and $\a$-dimensional elements remains
intact. Let $R$ denote the set of all permutations of $\a$. Define
$\Bb$ as the direct product of all the $\rho$-reducts of $\Aa$ for
$\rho\in R$, i.e.,
\[ \Bb \de \prod\langle \Rd^\rho\Aa : \rho\in R\rangle .\]
Clearly,  $\Bb$ is symmetric and representable. We show that it is
not inductive. Take the equation $e$ used in the proof of
Theorem~\ref{manyvar-thm}. We have seen in Remark~\ref{constr-rem}
that $\rho(e(c_0x))$ is valid in $\Aa$ for all $\rho\in R$. Hence,
$e(c_0x)$ is valid in all $\Rd^\rho\Aa$, hence in $\Bb$ by its
construction. However, $e(x)$ is not valid in $\Bb$ since it is not
valid in $\Aa$. This shows that $\Bb$ is not inductive. \QED

It is known that the same universal formulas are valid in $\Lf_\a$
as in $\RCA_\a$, see \cite[4.1.29]{HMTII}. There is no existential
formula distinguishing $\Lf_\a$ and $\RCA_\a$, either because each
$\RCA_\a$ has a subalgebra in $\Lf_\a$. The next complexity class is
$\Delta_2$-formulas, and our theorems so far imply that $\Lf_\a$
indeed can be distinguished from $\RCA_\a$ by a $\Delta_2$-formula.
We note that it was known that there is a $\Pi_2$-formula
distinguishing $\Lf_\a$ and $\RCA_\a$ (see \cite[2.6.53]{HMTI}).

\begin{cor}\label{delta-cor}
There is a $\Delta_2$-formula which is valid in $\Lf_\a$ but is not
valid in $\RCA_\a$.
\end{cor}

\noindent{\bf Proof.} The property of being inductive is defined by
a set $D$ of formulas of form $\forall\bar{x}e_1(\bar{x})\to
\forall\bar{x}e_2(\bar{x})$ where $e_1, e_2$ are equations using
variables occurring in $\bar{x}$. All such formulas are known to be
$\Delta_2$. Indeed, let $\varphi$ denote the previous formula. Then
$\varphi$ is equivalent both to the $\Pi_2$-formula
$\forall\bar{x}\exists\bar{y}(\neg e_2(\bar{x})\to \neg
e_1(\bar{y}))$, and to the $\Sigma_2$ formula
$\exists\bar{x}\forall\bar{y}(\neg e_1(\bar{x})\lor e_2(\bar{y}))$.
There is a representable algebra $\Aa$ which is not inductive, by
Theorem~\ref{ret-thm}(i). Since $\Aa$ is not inductive, there is a
$\Delta_2$-formula $\varphi\in D$ which is not valid in $\Aa$. Then
$\RCA_\a\not\models\varphi$ by $\Aa\in\RCA_\a$. However,
$\Lf_\a\models\varphi$ since $\Lf_\a\subseteq\Ret_\a$ by
Theorem~\ref{ret-thm}(ii).\QED

\begin{rem}\label{delta-rem} \rm{(i) We can get a concrete $\Delta_2$
formula separating $\Lf_\a$ and $\RCA_\a$ by using
Remark~\ref{constr-rem}(ii).

(ii) From the fact that there are more subvarieties of $\RCA_\a$
than generated by $\Lf_\a\subseteq\RCA_\a$ we can immediately get
that there is a subvariety  $\V$ of $\RCA_\a$ which is not generated
by its $\Lf_\a$-members, i.e., $\V$ is not generated by
$\V\cap\Lf_\a$.

(iii) Using (ii) above, from the fact that there are more
subvarieties of $\RCA_\a$ than generated by $\Lf_\a\subseteq\RCA_\a$
we can immediately get that there is a $\Delta_2$-formula
distinguishing $\Lf_\a$ and $\RCA_\a$, because the structure of
subvarieties of a variety $\V$ is determined by its
$\Delta_2$-theory. Indeed, assume that $\K$ and $\LL$ have the same
$\Delta_2$-theories. Then the same varieties are generated by
subclasses of $\K$ and $\LL$, since all of the formulas of form
$\forall\bar{x}e_1(\bar{x})\land\dots\land\forall\bar{x}e_n(\bar{x})\to
\forall\bar{x}e_0(\bar{x})$ are $\Delta_2$. Indeed, let
$\K_0\subseteq\K$, let $E_0=\Eq(\K_0)$ and let $\LL_0=\{\Aa\in\LL :
\Aa\models E_0\}$, then $\Eq(\LL_0)=E_0$, since
$E_0\subseteq\Eq(\LL_0)$ by the definition of $\LL_0$ and for all
$e\notin E_0$ we have $\K\not\models \Sigma\to e$ for all finite
$\Sigma\subseteq E$, so the same is true for $\LL$.

(iv) From what we said so far, it follows that for any
$\Lf_\a\subseteq\K\subseteq\Sy\cap\RCA_\a$ we have that the
$\Delta_2$-theories of $\K$ and $\RCA_\a$ are different but the
corresponding universal and existential theories coincide.}
\end{rem}

\subsection{Characterization of the equational theory of $\RCA_\a$}\label{rcaeq-ssec}

In this section we concentrate on sets of equations, rather than on
algebras. Assume $E$ is a set of equations in the language of
$\CA_\a$, it contains the cylindric axioms $(C_0) - (C_7)$
axiomatizing $\CA_\a$ and it is semantically closed (i.e., $e\in E$
iff $E\models e$). We call $E$ \emph{inductive} iff
$e(c_ix_1,\dots,c_ix_n)\in E$ implies $e(x_1,\dots,x_n)\in E$
whenever $i$ does not occur as an index of an operation in $e$.
Thus, an algebra is inductive iff its equational theory is such.
However, we will see that not all models of an inductive set of
equations are inductive. In the next theorem we characterize the
inductive sets of equations. We obtain that they coincide with the
equational theories of subclasses of $\Lf_\a$. Equational theories
of subclasses of $\Lf_\a$ are important, because $\Lf_\a$s
correspond to ordinary first order logic theories
(\cite[4.3.28(iii)]{HMTII}).

\begin{thm}\label{indset-thm}\hspace{10pt}
\[ E\text{ is inductive }\quad\text{ iff }\quad E=\Eq(\K)\text{ for
some }\K\subseteq\Lf_\a .\]
\end{thm}

\noindent{\bf Proof.} Assume that $E$ is inductive. Let $\Ff$ be the
$E$-free $\omega$-generated algebra. Then $E=\Eq(\Ff)$ and $\Ff$ is
inductive, by $E$ being inductive. So, there is $\Bb\in\Lf_\a$ with
$\Eq(\Ff)=\Eq(\Bb)$, by Theorem~\ref{ret-thm}(ii). This shows that
$E=\Eq(\K)$ for $\K = \{\Bb\}\subseteq\Lf_\a$. Assume now
$\K\subseteq\Lf_\a$ and let $E=\Eq(\K)$. Then $E$ contains the
cylindric axioms $(C_0)-(C_7)$ and is semantically closed. Also, $E$
is inductive by Theorem~\ref{ret-thm}(ii).\QED

Let us call \emph{inductive rule} the rule according to which from
$e(c_ix_1,\dots,c_ix_n)$ we can infer $e(x_1,\dots, x_n)$ provided
that $i\notin\ind(e)$. Note that this is a decidable rule, because
given any equation we can decide whether it is of form
$e(c_ix_1,\dots,c_ix_n)$ for an equation $e$ such that
$i\notin\ind(e)$.

\begin{cor}\label{rcaeq-cor}
The equational theory of $\RCA_\a$ is the least set of equations
which \\
contains the equations $(C_0)-(C_7)$ which define $\CA_\a$,\\
is closed under the 5 rules of equational logic,\ \ and\\
is closed under the inductive rule defined above.
\end{cor}

\noindent{\bf Proof.} By definition, a set $E$ of equations contains
$(C_0)-(C_7)$, is closed under the 5 rules of equational logic, and
is closed under the inductive rule iff $E$ is inductive. This is so
because equational logic is complete for its five rules. By
Theorem~\ref{indset-thm}, the least such set axiomatizes the variety
generated by the largest subclass of $\Lf_\a$, which subclass is
$\Lf_\a$ itself. Now, the variety generated by $\Lf_\a$ is
$\RCA_\a$, e.g., by \cite[4.1.29]{HMTII}. \QED

Corollary~\ref{rcaeq-cor} above gives a simple, natural enumeration
for the equational theory of $\RCA_\a$. It can be considered as a
solution to \cite[Problem 4.1]{HMTII} which asks for a simple
equational base for $\Eq(\RCA_\a)$. Certainly, the enumeration based
on the above Corollary~\ref{rcaeq-cor} is much simpler than any of
the three such enumerations given in \cite[sec.4.1]{HMTII}. It has
some resemblance to the second and third enumerations given in
\cite{HMTII}. An advantage of the present enumeration is that it
stays strictly in the equational language of $\CA_\a$ while the
second method given in \cite{HMTII} uses all first order logic
formulas in the language of $\CA_\a$, and the third method even uses
symbols outside the language of $\CA_\a$. A drawback of the present
enumeration is that it works only for infinite $\a$, while the three
methods given in \cite{HMTII} work for finite $\a$ also. We note
that possible solutions for Problems 4.1 and the related Problem
4.16 were also given in Simon~\cite{Simon} and Venema~\cite{Venema}.
The root of \cite[Problem 4.1]{HMTII} is Monk's theorem saying that
$\RCA_\a$ is not finite schema axiomatisable, exposing a gap between
abstract and representable cylindric algebras. As we mentioned in
the Introduction, this gap is addressed many ways in algebraic
logic, some works in this direction are \cite{AT, HH, HH97, SaGyu,
Sain, Ferenczi}.

\begin{rem}\label{singlelf-rem} \rm{
(i) Not all models of an inductive set of equations are inductive.
An example is $\Eq(\RCA_\a)$. It is inductive because $\RCA_\a$ is
generated by $\Lf_\a$ and it has a noninductive algebra by
Theorem~\ref{ret-thm}(i). Exceptions are the equational theories of
the minimal cylindric algebras in the sense that all members of
these varieties are inductive. We wonder whether these are the only
such exceptions or not.

(ii) Any variety of cylindric algebras generated by a class of
locally finite dimensional algebras is also generated by a single
$\Lf_\a$. This was known, but this also follows from
Theorems~\ref{ret-thm}, \ref{indset-thm} as follows. Let $\V$ be
generated by $\K\subseteq\Lf_\a$. Then $\Eq(\V)$ is inductive by
Theorem~\ref{indset-thm}, so the free algebra $\Ff$ of $\V$ is
inductive, then it is equationally indistinguishable from a
$\Bb\in\Lf_\a$ by Theorem~\ref{ret-thm}, and then
$\Eq(V)=\Eq(\Bb)$.

(iii) A set $E$ is inductive iff there is an ordinary first order
logic theory $\Th$ such that $E$ is the equational theory of all the
concept algebras of models of $\Th$. We briefly sketch a proof for
this, we deal with the logical connections in detail in another
paper. Let  $E$ be any inductive set. Then, by (ii) above, it is the
equational theory of a single $\Bb\in\Lf_\a$. Then $\Bb$ is the
Lindenbaum-Tarski algebra of an ordinary theory $\Th$, by
\cite[4.3.28(ii)]{HMTII}. It is not difficult to see that the
Lindenbaum-Tarski algebra is a subdirect product of $\{ \mbox{\sf
Ca}^{\Mm} : \Mm\models\Th\}$, which finishes the proof.}
%
\end{rem}

\subsection{Continuum many inductive varieties}\label{fewvar-ssec}

We close the paper with showing that subclasses of concept algebras
of ordinary first order logic with infinite universes generate
continuum many subvarieties. The proof of this theorem will be
analogous to, but simpler than, the proof of
Theorem~\ref{manyvar-thm}.
Concept algebras of ordinary first order logic with finite universes
also generate continuum many varieties, a slightly modified version
of the proof of \cite[4.1.24]{HMTII} shows this. This is why we deal
with concept algebras of models with infinite universes below.

Let $a_m\de c_{(m)}\prod\{-d_{ij} : i<j<m\}$, for $m\in\omega$, cf.\
\cite[2.4.61]{HMTI}. We call a cylindric algebra \emph{of infinite
base} iff $\{ e_m : m\in\omega\}$ is valid in it, and
${}_\infty\Lf_\a$ denotes the class of $\Lf_\a$s of infinite bases.
An \emph{inductive variety of infinite base} is a variety whose
equational theory is inductive and which contains the equations $\{
a_m=1 : m\in\omega\}$. The inductive varieties of infinite base are
exactly the varieties generated by subclasses of ${}_\infty\Lf_\a$,
by Theorem~\ref{ret-thm}. Also, they are exactly the varieties
generated by concept algebras of ordinary first order logic with
infinite bases, by Remark~\ref{singlelf-rem}(iii).

The following is a counterpoint to Theorem~\ref{manyvar-thm}. We
know that there can be only continuum many inductive varieties for
all $\a$ because inductive varieties are also symmetric. The
following theorem says that there are indeed continuum many
inductive varieties for all $\a$, even if we require the bases to be
infinite.

\begin{thm}\label{fewvar-thm}
Subclasses of ${}_{\infty}\Lf_a$ generate continuum many
subvarieties, for all infinite $\a$. In other words, there are
continuum many inductive varieties of infinite base.
\end{thm}

\noindent{\bf Proof.} As in the proof of Theorem~\ref{manyvar-thm},
we will use a set of independent equations, in this case we will use
a countable set of independent equations. The $n$-th equation $e_n$
will express that there is no partition of the universe (in the form
of an equivalence relation as element of the algebra) all of whose
blocks have size $n$. Then, for each $n\in\omega$ we will exhibit an
algebra $\Aa_n\in{}_{\infty}\Lf_\a$ in which $e_n$ fails, but $e_k$
holds for all $k\in\omega-2$, $k\ne n$.

We begin to write up the term expressing that ``$x$ is not an
equivalence relation on the whole base set with each equivalence
block having size $n$". The following terms express the parts of
this statement (in the final equation we will replace $x$ with
$c_2\dots c_{n}x$). Let $n\ge 2$. \smallskip

\noindent The domain of $x$ is not the base set:
\[ \delta(x) \de c_0-c_1x .\]
$x$ is not symmetric:
\[ \sigma(x) \de c_0c_1({}_2\s(0,1)x\oplus x) .\]
$x$ is not transitive:
\[ \tau(x) \de c_0c_1c_2(x\cdot \s^{01}_{12}x-\s^{01}_{02}x) .\]
$x$ is not reflexive:
\[ \rho(x) \de c_0c_1(d_{01}-x) .\]
There is a block in $x$ with size $<n$:
\[ \mu_<(x) \de c_0-c_1\dots -c_{n-1}(\prod\{-d_{ij} : i<j<n\}\cdot\prod\{ \s^{01}_{ij}x : i<j<n\}) .\]
There is a block in $x$ with size $>n$:
\[ \mu_>(x) \de c_{(n+1)}(\prod\{-d_{ij} : i<j\le n\}\cdot\prod\{ \s^{01}_{ij}x : i<j\le n\}) .\]
The sum of all these is
\[ \eta(x) \de \delta(x)+\sigma(x)+\tau(x)+\rho(x)+\mu_<(x)+\mu_>(x) .\]
The equation $e_n$ is defined as
\[ e_n(x)\quad \de\quad \eta(c_2\dots c_{n}x)=1 .\]

\begin{lem}\label{equiv-lem}
Let $\Aa\in\Cs_\a$, let $n\ge 2$ and let $a=c_2\dots c_na\in A$.
Then $\Aa\models e_n(a)$ iff for all $s\in a$ it is true that
$a[s,01]$ is not an equivalence relation on the base set with each
block of size $n$.
\end{lem}

\noindent{\bf Proof of Lemma~\ref{equiv-lem}.} Assume the conditions
of the lemma, then $\Aa\models e_n(a)$ iff for all $s\in{}^\a U$,
where $U$ is the base set of $\Aa$, we have
$s\in\eta(a)=\delta(a)+\dots + \mu_>(a)$. Let $R\de a[s,01]=\{ (u,v)
: s(01/uv)\in a\}\subseteq U\times U$. We have
\begin{equation}
 s\notin\delta(a)\quad\mbox{ iff }\quad\text{ the
domain of $R$ is $U$}.
\end{equation}
Indeed, $s\notin\delta(a)$ iff $s\in -\delta(a) = -c_0-c_1a$ iff for
all $u\in U$ there is $v\in U$ with  $s(01/uv)\in a$, which means
$(u,v)\in a[s,01]$.
\begin{equation}
 s\notin\sigma(a)\quad\mbox{ iff }\quad R\text{ is symmetric}.
\end{equation}
Indeed, $s\notin\sigma(a)=c_0c_1({}_2\s(0,1)a\oplus a)$ iff for all
$u,v\in U$ we have $s(01/uv)\notin({}_2\s(0,1)a\oplus a)$, this last
thing holds iff $s(01/vu)\in a\Leftrightarrow s(01/uv)\in a$, which
means that $R$ is symmetric.
\begin{equation}
 s\notin\tau(a)\quad\mbox{ iff }\quad R\text{ is transitive}.
\end{equation}
Indeed, $s\notin\tau(a) = c_0c_1c_2(a\cdot
\s^{01}_{12}a-\s^{01}_{02}a)$ iff for all $u,v,w\in U$ whenever
$s(012/uvw)\in a\cdot \s^{01}_{12}a$ we have $s(012/uvw)\in
\s^{01}_{02}a$. Now, $s(012/uvw)\in a\cdot \s^{01}_{12}a$ means that
$(u,v)\in R$ and $(v,w)\in R$ (we used $c_2a=a$). Similarly,
$s(012/uvw)\in \s^{01}_{02}a$ means that $(u,w)\in R$. Putting these
together, we get that $R$ is transitive.
\begin{equation}
 s\notin\rho(a)\quad\mbox{ iff }\quad R\text{ is reflexive}.
\end{equation}
Indeed, $s\notin\rho(a) = c_0c_1(d_{01}-a)$ iff for all $u,v\in U$
we have $u=v$ implies $(u,v)\in R$, i.e., $R$ is reflexive.

Assume now that $s\notin(\delta(a)+\sigma(a)+\tau(a)+\rho(a))$.
Then, by the above, we have that $R$ is an equivalence relation on
$U$.
\begin{equation}
 s\in\mu_<(a)\quad\mbox{ iff }\quad \text{ there is a block in $R$ with size }< n.
\end{equation}
Indeed, $s\in\mu_<(a) = c_0-c_1\dots -c_{n-1}(\prod\{-d_{ij} :
i<j<n\}\cdot\prod\{ \s^{01}_{ij}a : i<j<n\})$ iff there is $u_0\in
U$ such that there are no $u_1,\dots,u_{n-1}\in U$ such that
$u_0,\dots,u_{n-1}$ are all distinct and $(u_i,u_j)\in R$ for all
$i<j<n$. This means that there is a block in $R$ with size $<n$.
\begin{equation}
 s\in\mu_>(a)\quad\mbox{ iff }\quad \text{ there is a block in $R$ with size }> n.
\end{equation}
Indeed, $s\in\mu_>(a) = c_{(n+1)}(\prod\{-d_{ij} : i<j\le
n\}\cdot\prod\{ \s^{01}_{ij}a : i<j\le n\})$ iff there are
$u_0,\dots, u_n\in U$ such that they are all distinct and
$(u_i,u_j)\in R$ for all $i<j\le n$, and this means that there is a
block in $R$ with size $> n$.

By the above we have that $s\in\eta(a)$ iff whenever $R=s[a,01]$ is
an equivalence relation on $U$, there is either a block with size $<
n$ or else there is a block with size $>n$. This proves
Lemma~\ref{equiv-lem}.\bigskip

Let $n\in\omega, n\ge 2$, let $U$ be an infinite set, and let $R$ be
an equivalence relation on $U$ with each block of size $n$. Let
$\Aa_n$ be the $\Cs_\a$ with base $U$ and generated by $g\de\{ s\in
U^{\a} : (s_0,s_1)\in R\}$.

\begin{lem}\label{lfindep-lem}
$\Aa_n\not\models e_n$ but $\Aa_n\models e_k$ for all $k\ne n$,
$k\in\omega-2$.
\end{lem}

\noindent{\bf Proof of Lemma~\ref{lfindep-lem}.} $\Aa_n\not\models
e_n$ by Lemma~\ref{equiv-lem} and $g\in A_n$, since clearly
$s[g,01]=R$ for all $s\in g$ and $R$ is an equivalence relation of
the kind $e_n$ prohibits. Let $k\in\omega-2$, $k\ne n$, we want to
show that $\Aa_n\models e_k$. By Lemma~\ref{equiv-lem}, it is enough
to show that $a[s,01]$ is not an equivalence relation on $U$ with
all blocks of size $k$, whenever $s\in a = c_2\dots c_na$. We begin
doing this.

We call an $X\subseteq {}^\a U$ \emph{regular} if, intuitively, $X$
is determined by its restriction to its dimension set $\Delta(X)$,
formally
\[ s\in X\mbox{\ \ iff\ \ }z\in X,\quad \mbox{whenever $s,z\in{}^\a
U$ and $s,z$ agree on $\Delta(X)$}\] where $\Delta(X)\de\{ i\in\a :
c_iX\ne X\}$.
Since $\Aa_n$ is generated by $g$ which is a locally finite regular
element, we have that $a$ is also a regular locally finite element,
by \cite[3.1.64]{HMTII}. Let $S'\de \Rg(s\upharpoonright
\Delta(a)\cup\{ 0\})$ and let $S\de \{ u\in U : (\exists v\in
S')(u,v)\in R\}$. Then $S$ is finite since $\Delta(a)$ is finite and
each block of $R$ is finite. Assume that $E\de a[s,01]$ is an
equivalence relation on $U$ with all blocks finite, and $\ge 2$. If
$E$ is not such then we are done by Lemma~\ref{equiv-lem} and $k\ge
2$.

Let $u,v\in U-S$ such that $(u,v)\in E-R$, we will derive a
contradiction. Let $w\in U-S-u/R$ be arbitrary. There are infinitely
many such $w$. We want to show that $w\in u/E$, contradicting our
assumption that $u/E$ is finite. Let $\pi:U\to U$ be a permutation
of $U$ which leaves $R$ fixed, is identity on $S\cup \{ u\}$ and
takes $v$ to $w$. There is such a permutation since $v/R\cup w/R$ is
disjoint from $S\cup \{ u\}$ by our assumptions. Since $\pi$ leaves
$R$ fixed and $\Aa_n$ is generated by $g$, we have that $a$ is
closed under $\pi$, i.e., $z\in a$ iff $\pi\circ z\in a$ for all
$z$.  Now, $(u,v)\in E = a[s,01]$ means that $s(01/uv)\in a$.
Therefore $z\de \pi\circ (s(01/uv))\in a$. Now, $z_0=\pi(u)=u$,
$z_1=\pi(v)=w$, and $z$ agrees with $s(01/uw)$ on $\Delta(a)$ by
$\pi$ being the identity on $S$. Hence $s(01/uw)\in a$ by $z\in a$
and $a$ being regular. This means $(u,w)\in E$, i.e., $w\in u/R$ as
was to be shown.

Assume now that $u,v\in U-S$ such that $(u,v)\in R-E$, we will
derive a contradiction. Let $(w,v)\in E$, $w\ne v$. There is such by
our assumption that each block of $E$ has at least two elements. Let
$\pi: U\to U$ be a permutation of $U$ which is identity on $U-\{
u,v\}$ and interchanges $u$ and $v$. This $\pi$ is identity on $S$
by $u,v\notin S$ and it leaves $R$ fixed by $(u,v)\in R$. Thus, $a$
is closed under this $\pi$, too. As before, $(w,v)\in E = a[s,01]$
means that $s(01/wv)\in a$, therefore $z\de \pi\circ (s(01/wv))\in
a$. Then $z_0=\pi(w)=w$, $z_1=\pi(v)=u$, and $z$ agrees with
$s(01/wu)$ on $\Delta(a)$ by $\pi$ being the identity on $S$. Hence
$s(01/wu)\in a$ by $z\in a$ and $a$ being regular. This means
$(w,u)\in E$, contradicting $(u,v)\notin E$ and $(v,w)\in E$.

We have seen that $R$ and $E$ agree on the infinite set $U-S$. Since
each block of $R$ has $n$ elements, this means that $E$ has at least
one block with exactly $n$ elements. So, $e_k(a)$ holds in $\Aa$ by
$k\ne n$ and Lemma~\ref{equiv-lem}. By this, Lemma~\ref{lfindep-lem}
has been proved.\medskip

We are ready for completing the proof of Theorem~\ref{fewvar-thm}.
For each $H\subseteq\omega-2$ let $\V_H$ be the variety generated by
$\K_H\de\{\Aa_n : n\in H\}\subseteq{}_\infty\Lf_\a$. Assume
$G,H\subseteq\omega-2$ are distinct, say $n\in H-G$. Then
$\Aa_n\in\V_G-\V_H$ by Lemma~\ref{lfindep-lem}, so $\V_H$ and $\V_G$
are distinct. This shows that there are at least continuum many
varieties generated by subclasses of  ${}_\infty\Lf_\a$.  \QED

\bigskip\bigskip\bigskip

\noindent MTA Alfr\'ed R\'enyi Institute of Mathematics\\
Budapest, Re\'altanoda st.\ 13-15, H-1053 Hungary\\
andreka.hajnal@renyi.mta.hu, nemeti.istvan@renyi.mta.hu


\begin{thebibliography}{}

\bibitem{ACMNS} Andr\'eka, H., Comer, S. D., Madar\'asz, J. X., N\'emeti, I.,
and Sayed-Ahmed, T., Epimorphisms in cylindric algebras and
definability in finite variable logic. Algebra Universalis 61,3-4
(2009), 261-282.

\bibitem{AFerNCAIII} Andr\'eka, H., Ferenczi, M., and N\'emeti, I.,
(eds), Cylindric-like algebras and algebraic logic. Bolyai Society
Mathematical Studies  Vol 22, Springer Verlag Berlin, 2013. 478pp.

\bibitem{AGiNJSL94} Andr\'eka, H., Givant, S. R., and N\'emeti, I.,
The lattice of varieties of representable relation algebras. The
Journal of Symbolic Logic  59,2 (1994), 631-661.

\bibitem{AMonkN91} Andr\'eka, H., Monk, J. D., and N\'emeti, I., (eds),
Algebraic Logic. Colloq. Math. Soc. J. Bolyai Vol 54, North-Holland
Amsterdam,  1991. vi + 746pp.

\bibitem{ANComp} Andr\'eka, H., and N\'emeti, I., Comparing theories: the dynamics
of changing vocabulary. In: Johan V. A. K. van Benthem on logical
and informational dynamics. A. Baltag and S. Smets. eds., Springer
Series Outstanding contributions to logic Vol 5, Springer Verlag,
2014. pp.143-172.

\bibitem{AT} Andr\'eka, H., and Thompson, R. J., A Stone-type representation
theorem for algebras of relations of higher rank.  Trans. Amer.
Math. Soc. 309,2 (1988), 671-682.

\bibitem{ASSS13} Assem, M., Sayed-Ahmed, T., S\'agi, G., and Szir\'aki, D.,
The number of countable models via agebraic logic. Manuscript 2013.
http://real.mtak.hu/id/eprint/16985

\bibitem{BaHal15} Barrett, T. W., and Halvorson, H., Morita equivalence.
Preprint, arXiv:1506.04675, 2015.

\bibitem{BezhCAIII} Bezhanishvili, N., Varieties of two-dimensional cylindric algebras.
In: \cite{AFerNCAIII}, pp.37-60.

\bibitem{Blok} Blok, W. J., Varieties of interior algebras. PhD
Dissertation, University of Amsterdam, 1976.

\bibitem{Blok2} Blok, W. J., The lattice of modal logics: an
algebraic investigation. The Journal of Symbolic Logic 45 (1980),
221-236.


\bibitem{Craig} Craig, W. C., Logic in algebraic form: three
languages and theories. North-Holland Amsterdam, 1974. 204pp.

\bibitem{Ferenczi} Ferenczi, M., The polyadic generalization of the
Boolean axiomatization of fields of sets. Transactions of American
Mathematical Society 364,2 (2012), 867-896.

\bibitem{GKWZ} Gabbay, D. M., Kurucz, A., Wolter, F., and Zakharyaschev, M.,
Many-dimensional modal logics: theory and applications. Elsevier, 2003.

\bibitem{Givantbook} Givant, S. R., Relation algebras, vol I: Arithmetic and algebra;
vol II: Complete extensions, representations, varieties, and atom
structures. Springer-Verlag, to appear in 2016.

\bibitem{Goldblatt} Goldblatt, R., Varieties of complex algebras.
Annals of Pure and Applied Logic 44,3 (1989), 173-242.

\bibitem{Hal62} Halmos, P. R., Algebraic Logic. Chelsea Publ. Co.,
New York, 1962. 271pp.

\bibitem{HMTI} Henkin, L., Monk, J. D., and Tarski, A., Cylindric
Algebras. Part I. North-Holland Amsterdam, 1971, 1985.

\bibitem{HMTII} Henkin, L., Monk, J. D., and Tarski, A., Cylindric
Algebras. Part II. North-Holland Amsterdam, 1985.

\bibitem{HMTAN} Henkin, L., Monk, J. D., Tarski, A., Andr\'eka, H., and N\'emeti, I.,
Cylindric Set Algebras. Lecture Notes in Mathematics Vol 883,
Springer-Verlag Berlin 1981.  vi + 323pp.

\bibitem{HH} Hirsch, R., and Hodkinson, I., Step-by-step - building
representations in algebraic logic. The Journal of Symbolic Logic,
62 (1997), 225-279.

\bibitem{HH97} Hirsch, R., and Hodkinson, R., Axiomatising various
classes of relation and cylindric algebras. Logic J. of IGPL 5
(1997), 209-229.

\bibitem{HHbook} Hirsch, R., and Hodkinson, I., Relation algebras by
games. North-Holland Amsterdam, 2002. 730pp.

\bibitem{JiRo} Jipsen, P., and Rose, H., Varieties of lattices. Lecture
Notes in Mathematics Vol 1533, Springer Verlag Berlin, 1992.

\bibitem{Jonsson} J\'onsson, B., Varieties of relation algebras.
Algebra Universalis 15 (1982), 273-298.

\bibitem{LefSzek} Lefever, K., and Sz\'ekely, G., Interpretation of
special relativity in the language of Newtonian kinematics. Logic,
Relativity and Beyond, 2nd International Conference, August 9-13
2015, Budapest, talk in the Symposium of Equivalences of Theories
part.
http://www.renyi.hu/conferences/lrb15/slides/LRB15-Lefever--Szekely.pdf

\bibitem{MadDis} Madar\'asz, J. X., Logic and relativity (in the
light of definability theory). PhD Dissertation, ELTE Budapest, 2002.
xviii+367pp.

\bibitem{MadSzek} Madar\'asz, J. X., and Sz\'ekely, G., Comparing
relativistic and Newtonian dynamics in first-order logic. In: The
Vienna Circle in Hungary, Ver\"offentlichungen des Instituts Wiener
Kreis Vol 16, A. M\'at\'e, M. R\'edei, and F. Stadler, eds.,
Springer Vienna 2011. pp.155-179.

\bibitem{Marx} Marx, M., and Venema, Y., Multi-dimensional modal logic.
Kluwer Adacemic Publishers Dordrecht Boston and London, 1997.
xiii+239pp.

\bibitem{Monk69} Monk, J. D., On the lattice of equational classes
of one- and two-dimensional polyadic algebras. Notices Amer. Math.
Soc. vol 16 (1969), p.183.

\bibitem{Monk70} Monk, J. D., On equational classes of algebraic versions of logic I.
Mathematica Scandinavica 27 (1970), 53-71.

\bibitem{Madduxbook} Maddux, R. D., Relation Algebras. Elsevier Amsterdam,
2006. xxvi+731pp.

\bibitem{NVarprep86} N\'emeti, I., Varieties of cylindric algebras.
Preprint, Budapest, 1985.

\bibitem{NAPAL87} N\'emeti, I., On varieties of cylindric algebras with applications to
logic. Annals of Pure and Applied Logic 36 (1987) 235-277.

\bibitem{Nsurv} N\'emeti, I., Algebraizations of quantifier logics, an introductory
overview. Studia Logica Special issue on Algebraic Logic (W. J.
Blok, and D. Pigozzi eds) 50, 3-4 (1991), 485-569.

\bibitem{Rybakov} Rybakov, V. V., Admissibility of logical interence rules.
Elsevier Amsterdam, 1997.

\bibitem{Sain} Sain, I., On the search for a finitizable
algebraization of first order logic. Logic Journal of the IGPL 8,4
(2000), 495-589.

\bibitem{SaGyu} Sain, I., and Gyuris, V., Finite schematizable algebraic
logic. Logic Journal of the IGPL 5,5 (1997), 699-675.

\bibitem{ST} Sain, I. and Thompson, R. J., Strictly finite schema axiomatization
of quasi-polyadic algebras. In: \cite{AMonkN91}, pp.539-571.

\bibitem{SagiSziraki} S\'agi, G., and Szir\'aki, D., Some variants of Vaught's conjecture
from the perspective of algebraic logic. Logic Journal of the IGPL
20,6 (2012), 1064-1082.

\bibitem{Sayed} Sayed-Ahmed, T., Splitting methods in algebraic
logic: proving results on non-atom-canonicity, non-finite
axiomatizability and non-first-order definability for cylindric and
relation algebras. Preprint arXiv:1503.02189, 2015.

\bibitem{Sereny} Ser\'eny, Gy., Isomorphisms of finite cylindic set
algebras of characteristic zero. Notre Dame Journal of Formal Logic
34,2 (1993), 284-294.

\bibitem{Simon} Simon, A., Finite schema completeness for typeless logic and
representable cylindric algebras. In: \cite{AMonkN91}, pp.665-670.

\bibitem{Venema1} Venema, Y., Many-dimensional modal logic.
PhD Dissertation, Amsterdam, 1992. 177pp.

\bibitem{Venema} Venema, Y., Cylindric modal logic. The Journal of
Symbolic Logic 60,2 (1995), 591-623.

\bibitem{Weatherall} Weatherall, J. O., Are Newtonian gravitation and
geometrized Newtonian gravitation theoretically equivalent?
Preprint, arXiv:1411.5757, 2014.

\end{thebibliography}
\end{document}